\documentclass[a4paper,10pt]{article}
\usepackage[utf8x]{inputenc}
\usepackage{tracefnt,amsmath,tabu,array}
\usepackage{amssymb,graphicx,setspace,amsfonts,amsbsy}
\usepackage{pifont,latexsym,ifthen,amsthm,rotating,calc,textcase,booktabs,cancel,slashed}

\newtheorem{theorem}{Theorem}[section]
\newtheorem{lemma}[theorem]{Lemma}
\newtheorem{corollary}[theorem]{Corollary}
\newtheorem{proposition}[theorem]{Proposition}

\newtheorem{remark}[theorem]{Remark}
\newcommand{\filledbox}{\leavevmode
  \hbox to.77778em{%
  \hfil\vbox to.675em{\hrule width.6em height.6em}\hfil}}

\newcommand{\Rm}{{\mathbb R}}

\newcommand{\eps}{\varepsilon}

\begin{document}
%\doublespacing
\tabulinesep=1.0mm
% Enter full title and short title for running headers
\title{An inequality concerning Radon transform and non-radiative linear waves via a geometric method}
%This work is supported by National Natural Science Foundation of China Programs 11601374, 11771325}}

\author{Liang Li, Ruipeng Shen and Chenhui Wang\\
Centre for Applied Mathematics\\
Tianjin University\\
Tianjin, China
}

\maketitle

\begin{abstract}
 In this work we consider the operator
 \[
  (\mathbf{T} G) (x)= \int_{\mathbb{S}^2} G(x\cdot \omega, \omega) d\omega, \quad x\in \Rm^3, \; G\in L^2(\Rm\times \mathbb{S}^2).
 \]
 This is the adjoint operator of the Radon transform. We manage to give an optimal $L^6$ decay estimate of $\mathbf{T} G$ near the infinity by a geometric method, if the function $G$ is compactly supported. As an application we give decay estimate of non-radiative solutions to the 3D linear wave equation in the exterior region $\{(x,t)\in \Rm^3 \times \Rm: |x|>R+|t|\}$. This kind of decay estimate is a key element of the channel of energy method for wave equations. 
\end{abstract}

\section{Introduction} \label{sec:intro}
\subsection{Background and topics}
In this article we consider an operator 
\begin{equation} \label{the operator}
 (\mathbf{T} G) (x)= \int_{\mathbb{S}^2} G(x\cdot \omega, \omega) d\omega, \quad x\in \Rm^3, \; G\in L^2(\Rm\times \mathbb{S}^2). 
\end{equation}
This operator is exactly the adjoint of the Radon transform defined by (Here $dS$ is the usual measure of the plane $\omega \cdot x = s$.)
\[
 (\mathcal{R} f) (s,\omega) = \int_{\omega\cdot x = s} f(x) dS(x), \qquad (s,\omega) \in \Rm \times \mathbb{S}^2.
\]
Namely we always have $\langle f, \mathbf{T} G\rangle = \langle \mathcal{R} f, G\rangle $ for suitable functions $f(x)$ and $G(s,\omega)$. The angled brackets here are the corresponding pairing in the spaces $\Rm^3$ and $\Rm \times \mathbb{S}^2$, respectively. The application of the Radon transform includes partial differential equations, X-ray technology and radio astronomy. More details about the Radon transforms can be found in in Helgason \cite{radon1, radonbook} and Ludwig \cite{radon2}. In this work we are mainly interested in the application of the operator $\mathbf{T}$ on the wave equations. 
This operator helps solve the free waves, i.e. the solutions to homogenous linear wave equation $\partial_t^2 u - \Delta u = 0$, from their corresponding radiation fields. 

\paragraph{Radiation field} Let us first give a brief introduction of the radiation fields. The conception of radiation field dates back to 1960's, see Friedlander \cite{radiation1, radiation2}. Generally speaking, radiation fields discuss the asymptotic behaviours of free waves as time goes to infinity. The following version of statement can be found in Duyckaerts-Kenig-Merle \cite{dkm3}.

\begin{theorem}[Radiation field] \label{radiation}
Assume that $d\geq 3$ and let $u$ be a solution to the free wave equation $\partial_t^2 u - \Delta u = 0$ with initial data $(u_0,u_1) \in \dot{H}^1 \times L^2(\Rm^d)$. Then ($u_r$ is the derivative in the radial direction)
\[
 \lim_{t\rightarrow \pm \infty} \int_{\Rm^d} \left(|\nabla u(x,t)|^2 - |u_r(x,t)|^2 + \frac{|u(x,t)|^2}{|x|^2}\right) dx = 0
\]
 and there exist two functions $G_\pm \in L^2(\Rm \times \mathbb{S}^{d-1})$ so that
\begin{align*}
 \lim_{t\rightarrow \pm\infty} \int_0^\infty \int_{\mathbb{S}^{d-1}} \left|r^{\frac{d-1}{2}} \partial_t u(r\theta, t) - G_\pm (r\mp t, \theta)\right|^2 d\theta dr &= 0;\\
 \lim_{t\rightarrow \pm\infty} \int_0^\infty \int_{\mathbb{S}^{d-1}} \left|r^{\frac{d-1}{2}} \partial_r u(r\theta, t) \pm G_\pm (r\mp t, \theta)\right|^2 d\theta dr & = 0.
\end{align*}
In addition, the maps $(u_0,u_1) \rightarrow \sqrt{2} G_\pm$ are bijective isometries from $\dot{H}^1 \times L^2(\Rm^d)$ to $L^2 (\Rm \times \mathbb{S}^{d-1})$. 
\end{theorem}
\paragraph{Explicit formula} We call the functions $G_\pm$ radiation profiles in this work. They can be viewed as the ``initial data'' of free waves at the time $t = \pm \infty$. We may give an explicit formula for the one-to-one map from radiation fields $G_-(s,\omega)$ back to the initial data $(u_0,u_1)$ in dimension $3$:
\begin{align*}
  u_0(x) & = \frac{1}{2\pi} \int_{\mathbb{S}^{2}}  G_- \left(x \cdot \omega, \omega\right) d\omega;\\
 u_1(x) & = \frac{1}{2\pi} \int_{\mathbb{S}^{2}}  \partial_s G_- \left(x \cdot \omega, \omega\right)  d\omega.
\end{align*}
A similar formula has been known for many years, see Friedlander \cite{inverseradiation}. One may also refer to Li-Shen-Wei \cite{shenradiation} for an explicit formula for all dimensions $d\geq 2$. This map between initial data and radiation profiles can also be given in term of their Fourier transforms, as given in a recent work C\^{o}te-Laurent \cite{newradiation}. A formula of free waves in term of the radiation fields $G_-$ immediately follows by a time translation
\begin{equation} \label{from radiation to free wave}
 u(x,t) = \frac{1}{2\pi} \int_{\mathbb{S}^{2}}  G_- \left(x \cdot \omega +t , \omega\right) d\omega.
\end{equation}
We recall that the map from the radiation fields $G_-$ to initial data $(u_0,u_1)$ is an isometry from $L^2 (\Rm \times \mathbb{S}^2)$ to $\dot{H}^1 \times L^2 (\Rm^3)$. Thus the formula $u_0 = (1/2\pi) \mathbf{T} G_-$ implies that the operator $\mathbf{T}$ is a bounded linear operator from $L^2(\Rm\times \mathbb{S}^2)$ to $\dot{H}^1(\Rm^3)$. We may combine this with the Sobolev embedding $\dot{H}^1(\Rm^3) \hookrightarrow L^6(\Rm^3)$ and obtain that
$\mathbf{T}$ is also a bounded operator from $L^2(\Rm \times \mathbb{S}^2)$ to $L^6 (\Rm^3)$. 

\paragraph{Non-radiative solutions} In this work we are particularly interested in the case when $G$ is compactly supported ($b\in \Rm^+$)
\[
 \hbox{Supp} \; G \subseteq [-b,b]\times \mathbb{S}^2. 
\]
These radiation profiles correspond to the non-radiative solutions of the linear wave equation. More precisely, $G$ is a radiation profile with compact support as above, if and only if the corresponding free wave $u(x,t)$ given by \eqref{from radiation to free wave} satisfies (see Li-Shen-Wei \cite{shenradiation}, for example)
\begin{equation} \label{non-radiative condition}
 \lim_{t\rightarrow \pm \infty} \int_{|x|>b+|t|} |\nabla_{t,x} u(x,t)|^2 dx = 0. 
\end{equation}
These solutions are usually called non-radiative solutions, or more precisely, $b$-weakly non-radiative solutions. They play an important role in the channel of energy method, which becomes a powerful tool in the study of asymptotic behaviour of solutions in the past decade. Generally speaking, channel of energy method discusses the energy of solutions to the linear and/or non-linear wave equation in the exterior region $\{x: |x|>R+|t|\}$ for a constant $R$ as $t\rightarrow \pm \infty$. The basic theory of this method can be found in C\^{o}te-Kenig-Schlag \cite{channeleven}, Duyckaerts-Kenig-Merle \cite{tkm1, oddtool} and Kenig-Lawrie-Schlag \cite{channel5d}, for example. The application of channel of energy method includes proof of the soliton resolution conjecture for radial solutions to focusing, energy critical wave equation in all odd dimensions $d\geq 3$ by Duyckaerts-Kenig-Merle \cite{se, oddhigh} and the non-existence of soliton-like minimal blow-up solution in the energy super-critical or sub-critical case by Duyckaerts-Kenig-Merle \cite{dkm2} and Shen \cite{shen2}, for instance. 
\paragraph{Decay estimate} One important part of channel of energy theory is to show that if $u$ is a non-radiative solution to a suitable non-linear wave equation, then the asymptotic behaviour of its initial data as $x\rightarrow +\infty$ is similar to that of non-radiative free waves. (see Duyckaerts-Kenig-Merle \cite{oddtool}, for example) The idea is to show that the nonlinear term gradually becomes negligible in the exterior region $\{(x,t)\in \Rm^3\times \Rm: |x|>R+|t|\}$ as $R\rightarrow +\infty$. As a result, this argument depends on suitable decay estimates of linear non-radiative free waves in the exterior region $\{(x,t): |x|>|t|+R\}$. Most previously known results of this kind depends on the radial assumption on the solutions. This work is an attempt to give a decay estimate as mentioned above in the non-radial case. This decay estimate is used in an accompanying paper to give the asymptotic behaviour of weakly non-radiative solutions to a wide range of non-linear wave equations, without the radial assumption. 

\paragraph{Topics} The main topic of this work is to find a good upper bound of the integral 
\[
 \int_{|x|>R} |\mathbf{T}G (x)|^6 dx 
\]
when the radiation profile $G$ is compactly supported. This immediately gives a decay estimate of non-radiative linear waves.  

\begin{remark}
 Strictly speaking, Li-Shen-Wei \cite{shenradiation} only gives proof of \eqref{from radiation to free wave} for smooth and compactly supported radiation fields $G_-$. But the same formula holds for any radiation profiles $G_- \in L^2(\Rm\times \mathbb{S}^2)$. More precisely, given any time $t$, the integral
 \[
  u(x,t) = \frac{1}{2\pi} \int_{\mathbb{S}^{2}}  G_- \left(x \cdot \omega +t , \omega\right) d\omega
 \]
 is defined for almost everywhere $x\in \Rm^3$ so that $u(x,t)$ is a linear free wave with radiation field $G_-$. In order to prove this we only need to use the result for smooth and compactly supported radiation fields and apply the classic approximation techniques of real analysis. 
\end{remark}

\subsection{Main results}

Now we give the statement of our main results. 

\begin{proposition} \label{main 1}
The linear operator $\mathbf{T}$ defined in \eqref{the operator} satisfies 
\begin{itemize}
 \item[(a)] Assume $b>a>0$ with $b/a\leq 2$. If $G\in L^2(\Rm\times \mathbb{S}^2)$ is supported in $([-b,-a]\cup [a,b])\times \mathbb{S}^2$, then we have
 \begin{align*}
  \int_{|x|>R} |\mathbf{T} G (x)|^6 dx & \lesssim  \frac{(a/R)^2 (1-a/b)^3}{1-a/R}  \|G\|_{L^2(\Rm\times \mathbb{S}^2)}^6, \qquad \forall R\geq b;\\
  \int_{\Rm^3} |\mathbf{T} G (x)|^6 dx & \lesssim (1-a/b)^2 \|G\|_{L^2(\Rm\times \mathbb{S}^2)}^6. 
 \end{align*}
\item[(b)] Assume $R \geq b>0$. If $G\in L^2(\Rm\times \mathbb{S}^2)$ is supported in $([-b,b])\times \mathbb{S}^2$, then
 \[
  \int_{|x|>R} |\mathbf{T} G (x)|^6 dx \lesssim (b/R)^2 \|G\|_{L^2(\Rm\times \mathbb{S}^2)}^6.
 \]
\end{itemize} 
\end{proposition}

\noindent We may utilize Proposition \ref{main 1} and conduct  a detailed discussion of various choices of parameters $a,b$ to obtain the following decay estimate
\begin{corollary} \label{cor1}
 If $G\in L^2(\Rm\times \mathbb{S}^2)$ is supported in $[a,b]\times \mathbb{S}^2$, then 
 \[
  \int_{|x|>R} |\mathbf{T} G (x)|^6 dx \lesssim  \frac{(b-a)^2}{R^2}  \|G\|_{L^2(\Rm\times \mathbb{S}^2)}^6.
 \] 
\end{corollary}

\begin{remark}
 The decay rate $\|\mathbf{T} G\|_{L^6(\{x: |x|>R\})} \lesssim R^{-1/3} \|G\|_{L^2}$ given above is optimal. We define ($\omega = (\omega^1, \omega^2, \omega^3) \in \mathbb{S}^2$)
\[
 G(s,\omega) = \left\{\begin{array}{ll} 1, & \hbox{if}\; s\in [-1,1], \; 0<\omega^3<\frac{1}{6R};\\
 0, & \hbox{otherwise;}  \end{array}\right.
\] 
A basic calculation shows that $\|G\|_{L^2} \simeq R^{-1/2}$. In addition, we have $|x\cdot \omega|\leq 1$ if $|x_3|<2R$, $x_1^2 + x_2^2 < 1/9$ and $\omega^3 \in (0,1/6R)$. Therefore  
\[
 \mathbf{T}G (x) = \int_{\omega \in \mathbb{S}^2, 0<\omega^3 < 1/6R} G(x\cdot \omega,\omega) d\omega \simeq R^{-1}, \qquad |x_3|<2R, \; x_1^2 + x_2^2 < 1/9. 
\]
Therefore we have $\|\mathbf{T} G\|_{L^6(\{x: |x|>R\})} \gtrsim R^{-5/6} \simeq R^{-1/3} \|G\|_{L^2}$. 
\end{remark}

\noindent The $L^6$ decay estimates of $\mathbf{T} G$ given above can be used to give decay estimate of non-radiative solutions in the exterior region

\begin{proposition} \label{main 2} 
 Let $u(x,t)$ be a solution to the 3-dimensional linear wave equation $\partial_t^2 u -\Delta u =0$ with a finite energy $E$ so that
 \[
  \lim_{t\rightarrow \pm \infty} \int_{|x|>r+|t|} |\nabla_{t,x} u(x,t)|^2 dx = 0. 
 \] 
Then the following inequalities hold
\begin{align*}
 &\|u\|_{L_t^\infty L^6 (\Rm \times \{x: |x|>R\})} \lesssim (r/R)^{1/3} E^{1/2},\; R\geq r;& &\|u(\cdot,t)\|_{L^6(\Rm^3)} \lesssim (r/|t|)^{1/3} E^{1/2}.&
\end{align*}
In addition, we have the decay of Strichartz norm in the exterior region
\[
 \|u\|_{L_t^p L^q (\Rm \times \{x:|x|>R\})} \lesssim_{p,q,\kappa} (r/R)^\kappa E^{1/2}, \qquad R \geq r; 
\]
as long as the constants $p, q \in (2,+\infty)$ and $\kappa >0$ satisfy $1/p+3/q=1/2$ and $\kappa < 2/q$.
\end{proposition}

\subsection{The idea}

The proof of our main result, Proposition \ref{main 1}, consists of four steps. 
\paragraph{Step 1} We temporarily assume $\hbox{Supp}\, G \subset [a,b]\times \mathbb{S}^2$ with $b>a>0$ and $b/a<2$. Other cases are direct consequences. We recall
\[
  u(x,t) = \int_{\mathbb{S}^2} G(x\cdot \omega+t, \omega) d\omega \in C(\Rm_t; \dot{H}^1(\Rm^3)) \hookrightarrow C(\Rm_t; L^6(\Rm^3)).
\]
Thus we may rewrite 
\[
 \mathbf{T} G = \lim_{\delta\rightarrow 0^+} \frac{1}{\delta} \int_{0}^\delta \int_{\mathbb{S}^2} G(x\cdot \omega+t,\omega) d\omega dt,
\]
and consider the upper bound of
\[
 \int_{|x|>R} \left|\frac{1}{\delta} \int_{0}^\delta \int_{\mathbb{S}^2} G(x\cdot \omega+t,\omega) d\omega dt\right|^6 dx. 
\]
A careful calculation gives an upper bound 
\begin{equation} \label{upper bound form 1} 
 K_\delta = \frac{1}{\delta^6} \int_{(\mathbb{S}^2 \times I)^6} \left(\prod_{k=1}^6 |G(s_k, \omega_k)| \right) |A_1 \cap A_2 \cap \cdots \cap A_6| (d\omega ds)^6 .
\end{equation}
Here $I = [a,b]$, $A_k = \{x\in \Rm^3: s_k-\delta < x\cdot \omega_k<s_k, |x|>R\}$ and 
\[
 (d\omega ds)^6 = \prod_{k=1}^6 d\omega_k ds_k.
\]
\paragraph{Step 2} In order to prove Proposition \ref{main 1}, we need to show $K_\delta \leq C\|G\|_{L^2}^6$ and determine the best constant $C$. The right hand side $\|G\|_{L^2}^6$ can be rewritten in the form 
\[
 \|G\|_{L^2(\Rm \times \mathbb{S}^2)}^6 = \int_{(\mathbb{S}^2 \times I)^3} |G(s_1,\omega_1)|^2 |G(s_2,\omega_2)|^2 |G(s_3,\omega_3)|^2 ds_{123} d\omega_{123}. 
\]
Here $ds_{ijk} = ds_i ds_j ds_k$, $d\omega_{ijk} = d\omega_i d\omega_j d\omega_k$. A comparison of this identity with \eqref{upper bound form 1} indicates that a Cauchy-Schwartz inequality might do the job. One could try to write (we define $A_{ijk} = A_i \cap A_j \cap A_k$ and $A_{123456} = A_1 \cap \cdots \cap A_6$)
\begin{align*}
 K_\delta & \leq \frac{1}{2\delta^6} \int_{(\mathbb{S}^2 \times I)^6} \frac{|A_{456}|\cdot |A_{123456}|}{|A_{123}|} |G(s_1,\omega_1)|^2  |G(s_2,\omega_2)|^2  |G(s_3,\omega_3)|^2  (ds d\omega)^6 \\
 & \qquad + \frac{1}{2\delta^6} \int_{(\mathbb{S}^2 \times I)^6} \frac{|A_{123}|\cdot |A_{123456}|}{|A_{456}|} |G(s_4,\omega_4)|^2  |G(s_5,\omega_5)|^2  |G(s_6,\omega_6)|^2  (ds d\omega)^6 \\
 & \leq \frac{1}{\delta^6} \int_{(\mathbb{S}^2 \times I)^6} \frac{|A_{456}|\cdot |A_{123456}|}{|A_{123}|} |G(s_1,\omega_1)|^2  |G(s_2,\omega_2)|^2  |G(s_3,\omega_3)|^2  (ds d\omega)^6\\
 & \leq \int_{(\mathbb{S}^2 \times I)^3} J(s_{123}, \omega_{123})|G(s_1,\omega_1)|^2 |G(s_2,\omega_2)|^2 |G(s_3,\omega_3)|^2 ds_{123} d\omega_{123};
\end{align*}
with 
\[
 J(s_{123}, \omega_{123}) = \frac{1}{\delta^6 |A_{123}|} \int_{(\mathbb{S}^2 \times I)^3} |A_{456}| \cdot |A_{123456}| d\omega_{456} ds_{456}.
\]
Here we put the weights $|A_{456}|/|A_{123}|$ for the purpose of balance, because without the weights the coefficient $J(s_{123},\omega_{123})$ would become
\[
 \frac{1}{\delta^6}\int_{(\mathbb{S}^2 \times I)^3} |A_{123456}| d s_{456} d\omega_{456},
\]
which seems to be proportional to $|A_{123}|$. Now we need to find an upper bound of 
\begin{align*}
 \sup_{\omega_{123}, s_{123}} J(s_{123}, \omega_{123})  = \sup_{\omega_{123}, s_{123}} \frac{1}{|A_{123}|} \int_{A_{123}} \left(\frac{1}{\delta^6}\int_{(\mathbb{S}^2 \times I)^3} |A_{456}| \chi_{A_{456}} (x) d\omega_{456} ds_{456}\right) dx.
\end{align*}
A reasonable upper bound can be found
\begin{align*}
 \sup_{x} \frac{1}{\delta^6}\int_{(\mathbb{S}^2 \times I)^3} |A_{456}| \chi_{A_{456}} (x) d\omega_{456} ds_{456} & = \sup_{x} \frac{1}{\delta^6}\int_{(\mathbb{S}^2 \times I)^3, x\in A_{456}} |A_{456}| d\omega_{456} ds_{456} \\
 & \leq \sup_{x} \int_{\Omega^3 (x)} \frac{1}{|(\omega_4\times \omega_5)\cdot \omega_6|} d\omega_{456}.
\end{align*}
Here $\Omega(x) = \{\omega\in \mathbb{S}^2: a-\delta < x\cdot \omega <b\}$. We use the following facts in the inequality above
\begin{itemize}
 \item $x\in A_k \Rightarrow \omega_k \in \Omega(x)$;
 \item $A_{456} \leq \delta^3 /|(\omega_4 \times \omega_5)\cdot \omega_6|$;
 \item Given $x$ and $\omega_k$, then $\{s_k \in [a,b]: x\in A_k\}$ is an interval whose length is no more than $\delta$. 
\end{itemize} 
Finally we need to find an upper bound of 
\[
 \sup_{x} \int_{\Omega^3 (x)} \frac{1}{|(\omega_4\times \omega_5)\cdot \omega_6|} d\omega_{456}.
\]
Unfortunately we have
\[
 \int_\Omega \frac{1}{|(\omega_4\times \omega_5)\cdot \omega_6|} d\omega_{456} = + \infty
\]
for any open region $\Omega \subset \mathbb{S}^2$. As a result, the argument above has to be improved in some way. The key observation here is that we have many different ways to split the product of $G(s_k,\omega_k)$ into two triples when we apply the Cauchy-Schwartz. In order to avoid too small value of $|(\omega_i \times \omega_j)\cdot \omega_k|$, which appears in the denominator in the integral given above, given $\omega_1, \omega_2, \cdots, \omega_6 \in \mathbb{S}^2$, we split them into two group of three, i.e. $(\omega_{k_1}, \omega_{k_2}, \omega_{k_3})$ and $(\omega_{k_4}, \omega_{k_5}, \omega_{k_6})$, so that the product
\[
 |(\omega_{k_1} \times \omega_{k_2})\cdot \omega_{k_3}|\cdot |(\omega_{k_4} \times \omega_{k_5})\cdot \omega_{k_6}|
\] 
gets its maximum value among all possible grouping method. In this work we call these kind of triples reciprocal triples. Following a similar argument as above but using reciprocal triples instead in the Cauchy Schwartz 
\begin{align*}
 \left(\prod_{k=1}^6 |G(s_k, \omega_k)| \right) & \leq \frac{|A_{k_4 k_5 k_6}|}{2|A_{k_1 k_2 k_3}|} |G(s_{k_1}, \omega_{k_1})|^2 |G(s_{k_2}, \omega_{k_2})|^2 |G(s_{k_3}, \omega_{k_3})|^2 \\
 & \qquad +  \frac{|A_{k_1 k_2 k_3}|}{2|A_{k_4 k_5 k_6}|} |G(s_{k_4}, \omega_{k_4})|^2 |G(s_{k_5}, \omega_{k_5})|^2 |G(s_{k_6}, \omega_{k_6})|^2;\\
 \{k_1,k_2, \cdots, k_6\} = \{1,2,\cdots, 6\}; & \quad (\omega_{k_1}, \omega_{k_2}, \omega_{k_3}), \;(\omega_{k_4}, \omega_{k_5}, \omega_{k_6})\; \hbox{are reciprocal},
\end{align*}
we reduce the problem to find an upper bound of 
\begin{equation} \label{least upper bound intro}
 \sup_{x\in B^c, \omega_{123} \in \Omega^3(x)}  \int_{\Sigma(\omega_{123})\cap \Omega^3(x)} \frac{1}{|(\omega_4\times\omega_5)\cdot \omega_6|} d\omega_{456}.
\end{equation}
Here $B^c=\{x: |x|>R\}$ and $\Sigma(\omega_{123}) \subseteq (\mathbb{S}^2)^3$ consists of all reciprocal triples of $\omega_{123}$. The reciprocal condition above significantly restricts the location, size and/or shape of the surface triangles $(\omega_4, \omega_5, \omega_6)$ thus leads to a finite least upper bound. The remaining work is to figure out this least upper bound. 

\paragraph{Step 3} We then apply a central projection $\mathbf{P} : \mathbb{S}_+^2 \rightarrow \Rm^2$ defined by $\mathbf{P} (x_1,x_2,x_3) = (x_1/x_3, x_2/x_3)$ and rewrite the least upper bound \eqref{least upper bound intro} in the form of an integral in Euclidean space $\Rm^2$: 
\[
 \sup_{x\in B^c} \left(\frac{b^6}{|x|^6} \sup_{Y_1,Y_2,Y_3 \in \Omega^\ast} \int_{\Sigma(Y_{123})\cap (\Omega^\ast)^3} \frac{1}{|\triangle Y_4 Y_5 Y_6|} dY_{456}\right).
\]
Here $\Omega^\ast$ is an annulus region (depending on $|x|$) in $\Rm^2$ and $\Sigma(Y_1 Y_2 Y_3)$ is the subset of $(\Rm^2)^3$ consisting of all reciprocal triples (or triangles) $Y_4 Y_5 Y_6$ in $\Rm^2$. Here reciprocal triangles in $\Rm^2$ are defined in a similar way to reciprocal triples in $\mathbb{S}^2$. 
\[
 |\triangle Y_1 Y_2 Y_3| \cdot |\triangle Y_4 Y_5 Y_6| \geq \frac{1}{65} \max_{\{k_1, k_2, \cdots, k_6\} = \{1,2,\cdots,6\}} |\triangle Y_{k_1} Y_{k_2} Y_{k_3}| \cdot |\triangle Y_{k_4} Y_{k_5} Y_{k_6}|.
\]

\paragraph{Step 4} In the final step we utilize the geometric properties of reciprocal triangles and give an upper bound
\begin{equation}\label{geometric inequality first}
 \sup_{Y_1, Y_2, Y_3 \in \Omega^\ast}  \int_{\Sigma(Y_{123})\cap (\Omega^\ast)^3} \frac{1}{|\triangle Y_4 Y_5 Y_6|} dY_{456} \lesssim w^3 r.
\end{equation} 
Here $r$ is the radius of outer boundary and $w$ is the width of the annulus region $\Omega^\ast$. Finally we may plug this upper bound in and conplete the proof of Proposition \ref{main 1}.

\subsection{Notations and Structure of this work}

\paragraph{Notations} In this work the notation $A \lesssim B$ means that there exists a constant $c$ so that $A \leq cB$. In this work these explicit constants are absolute constants, i.e. depends on nothing, unless stated otherwise. The notation $\gtrsim$ is similar. The meaning of $\ll$ is similar to $\lesssim$, i.e. there exists a constant $c$, so that $A \leq cB$. But in this case we additionally assume $c<1$ is very small. The meaning of $\gg$ is similar. We may add subscripts to these notations to indicate that the explicit constants depend on these subscripts but nothing else. Throughout this work we use the notation $\chi$ for characteristic functions and $|\Omega|$ for the Lebesgue measure of a subset $\Omega$ of the Euclidean spaces or the sphere $\mathbb{S}^2$.     
 
\paragraph{Structure of this work} In Section 2 we first reduce the proof of Proposition \ref{main 1} to a geometric inequality. Section 3 is devoted to the proof of some basic geometric properties regarding reciprocal triangles and circular annulus regions, which are the preparation work for the proof of the geometric inequality \eqref{geometric inequality first}. Next in Section 4 we prove the geometric inequality by considering reciprocal triangles with different sizes and angles separately. In Section 5 we combine all results from previous sections to finish the proof of Proposition \ref{main 1} and then give an application on the decay estimate of non-radiative solutions. Finally we give an estimate of Radon transform in the two or three dimensional space, as another application of our main result. 

\section{Transformation to a Geometric Inequality} \label{sec: transformation}

\begin{figure}[h]
 \centering
 \includegraphics[scale=1.25]{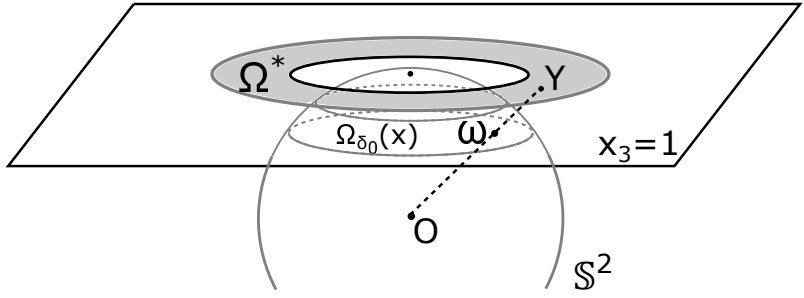}
 \caption{Illustration of projection $\mathbf{P}$} \label{figure central projection}
\end{figure}

\begin{figure}[h]
 \centering
 \includegraphics[scale=1.25]{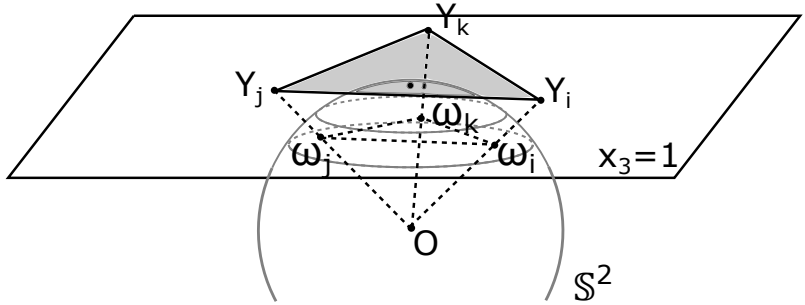}
 \caption{Illustration of volume} \label{figure projection2}
\end{figure}

In this section we reduce the proof of Proposition \ref{main 1} to a geometric inequality. Let us temporarily assume $G(s,\omega)$ is supported in $[a,b] \times \mathbb{S}^2$. Here $a,b>0$ so that $b/a\leq 2$. We recall that the function defined by 
\[
 u(x,t) = \int_{\mathbb{S}^2} G(x\cdot \omega+t, \omega) d\omega
\]
is a finite-energy free wave, thus we have 
\[
 u(\cdot,t) \in C (\Rm, \dot{H}^1(\Rm^3)) \quad \Rightarrow \quad u(\cdot,t) \in C (\Rm, L^6(\Rm^3)). 
\]
This immediately gives the following convergence in $L^6(\Rm^3)$
\[
 \lim_{\delta \rightarrow 0^+} \frac{1}{\delta} \int_0^\delta u(x,t) dt  = \mathbf{T} G.
\]
Thus it suffices to find an upper bound of 
\[
 \liminf_{\delta \rightarrow 0^+} \left\|\frac{1}{\delta} \int_0^\delta u(x,t) dt\right\|_{L^6(\{x: |x|>R\})}.
\]
We may rewrite 
\begin{align*}
 \frac{1}{\delta} \int_0^\delta u(x,t) dt & = \frac{1}{\delta} \int_0^\delta \int_{\mathbb{S}^2} G(x\cdot \omega +t, \omega) d\omega dt\\
 & = \frac{1}{\delta} \int_\Rm \int_{\mathbb{S}^2} G(s,\omega) \chi_{(0,\delta)} (s-x\cdot \omega) d\omega ds\\
 & = \frac{1}{\delta} \int_I \int_{\mathbb{S}^2} G(s,\omega) \chi_{(0,\delta)} (s-x\cdot \omega) d\omega ds
\end{align*}
Here $I =[a,b]$ and we use the compact-supported assumption of $G$. Given $\delta, s,\omega$, we may interpret $\chi_{(0,\delta)}(s-x\cdot \omega)$ as the characteristic function of the set 
\[
  A_{s,\omega,\delta} = \{x\in \Rm^3: s-\delta < x\cdot \omega<s\}.
\]
The set $A_{s, \omega, \delta}$ is a thin slice of the space $\Rm^3$, which is orthogonal to $\omega$ and a distance of about $s$ away from the origin. For convenience we introduce the notation $\chi_{s,\omega, \delta} (x) = \chi_{(0,\delta)}(s-x\cdot \omega)$. Thus we may rewrite 
\[
 \frac{1}{\delta} \int_0^\delta u(x,t) dt = \frac{1}{\delta} \int_I \int_{\mathbb{S}^2} G(s,\omega) \chi_{s, \omega, \delta} (x) d\omega ds
\]
Now we consider the integral
\[
  J_\delta = \int_{|x|>R} \left|\frac{1}{\delta} \int_0^\delta u(x,t) dt\right|^6 dx.
\]
We plug the explicit expression of $u$ in and obtain
\begin{align*}
 J_\delta & = \frac{1}{\delta^6} \int_{|x|>R} \left| \int_I \int_{\mathbb{S}^2} G(s,\omega) \chi_{s, \omega, \delta} (x) d\omega ds\right|^6 dx \\
 & \leq \frac{1}{\delta^6} \int_{|x|>R} \int_{(I\times \mathbb{S}^2)^6} \left(\prod_{k=1}^6 |G(s_k, \omega_k)| \chi_{s_k, \omega_k, \delta} (x) \right) (d\omega ds)^6 dx.
 \end{align*}
Here we slightly abuse the notation 
\[
 (d\omega ds)^6 = \prod_{k=1}^6 d\omega_k ds_k.
\]
Now we introduce reciprocal triples. If triples $(\omega_1, \omega_2, \omega_3), (\omega_4, \omega_5, \omega_6) \in (\mathbb{S}^2)^3$ satisfy
\[
 \left|(\omega_1 \times \omega_2) \cdot \omega_3\right| \left|(\omega_4\times \omega_5) \cdot \omega_6\right| = \max_{j_1, j_2, \cdots, j_6}  \left|(\omega_{j_1} \times \omega_{j_2}) \cdot \omega_{j_3}\right| \left|(\omega_{j_4}\times \omega_{j_5}) \cdot \omega_{j_6}\right|,
\]
we call these triples reciprocal to each other. Here the maximum is taken for all possible permutation of $\{1,2,\cdots,6\}$. By rotating the variables we only need to consider the integral in the region where the triples $(\omega_1, \omega_2, \omega_3), (\omega_4, \omega_5, \omega_6) \in (\mathbb{S}^2)^3$ are reciprocal. More precisely we have
\begin{align*}
 J_\delta \leq \frac{10}{\delta^6} \int_{|x|>R} \int_{\Sigma \times I^6} \left(\prod_{k=1}^6 |G(s_k, \omega_k)| \chi_{s_k, \omega_k, \delta} (x) \right) (d\omega ds)^6 dx.
\end{align*}
 Here 
\[
 \Sigma = \{(\omega_1, \cdots , \omega_6) \in (\mathbb{S}^2)^6:(\omega_1, \omega_2, \omega_3), (\omega_4, \omega_5, \omega_6)\; \hbox{are reciprocal}\}.
\]
For convenience we use the notations $A_k = A_{s_k, \omega_k, \delta} \cap \{x\in \Rm^3 : |x| >R\}$, $A_{ijk} = A_i \cap A_j \cap A_k$ and $A_{123456} = A_1 \cap \cdots \cap A_6$ below. We may rewrite 
\begin{align*}
 J_\delta \leq \frac{10}{\delta^6} \int_{|x|>R} \int_{\Sigma \times I^6} \left(\prod_{k=1}^6 |G(s_k, \omega_k)| \right) \chi_{A_{123456}}(x) (d\omega ds)^6 dx.
\end{align*}
We then apply Cauchy-Schwartz inequality 
\begin{align*}
 J_\delta & \leq \frac{5}{\delta^6} \int_{|x|>R} \int_{\Sigma \times I^6} \frac{|A_{456}|}{|A_{123}|} |G(s_1,\omega_1)|^2  |G(s_2,\omega_2)|^2  |G(s_3,\omega_3)|^2 \chi_{A_{123456}}(x) (d\omega ds)^6 dx\\
 & \qquad + \frac{5}{\delta^6} \int_{|x|>R} \int_{\Sigma \times I^6} \frac{|A_{123}|}{|A_{456}|} |G(s_4,\omega_4)|^2  |G(s_5,\omega_5)|^2  |G(s_6,\omega_6)|^2 \chi_{A_{123456}}(x) (d\omega ds)^6 dx\\
 & \leq \frac{10}{\delta^6} \int_{|x|>R} \int_{\Sigma \times I^6} \frac{|A_{456}|}{|A_{123}|} |G(s_1,\omega_1)|^2  |G(s_2,\omega_2)|^2  |G(s_3,\omega_3)|^2 \chi_{A_{123456}}(x) (d\omega ds)^6 dx.
\end{align*}
Next we use notations $d \omega_{ijk} = d\omega_i d\omega_j d\omega_k$, $ds_{ijk} = ds_i ds_j ds_k$ and rewrite the integral 
\begin{equation} \label{upper bound J delta}
 J_\delta \leq \frac{10}{\delta^6} \int_{(\mathbb{S}^2)^3 \times I^3} J(s_{123}, \omega_{123}) |G(s_1,\omega_1)|^2  |G(s_2,\omega_2)|^2  |G(s_3,\omega_3)|^2 ds_{123} d\omega_{123}.
\end{equation}
Here 
\begin{align*}
 J(s_{123}, \omega_{123}) & = \int_{|x|>R} \int_{\Sigma(\omega_{123})\times I^3} \frac{|A_{456}|}{|A_{123}|} \chi_{A_{123456}} (x) ds_{456} d\omega_{456} dx; \\
 \Sigma(\omega_{123}) & = \{(\omega_4, \omega_5, \omega_6)\in (\mathbb{S}^2)^3: (\omega_1, \omega_2, \omega_3), (\omega_4, \omega_5, \omega_6)\; \hbox{are reciprocal}\}.
\end{align*}
We may further find an upper bound of $J(s_{123}, \omega_{123})$. 
\begin{align*}
 J(s_{123}, \omega_{123}) 
 & = \frac{1}{|A_{123}|} \int_{A_{123}} \left(\int_{\Sigma(\omega_{123})\times I^3} |A_{456}| \chi_{A_{456}}(x) ds_{456} d\omega_{456}\right) dx \\
 & \leq \sup_{x\in A_{123}} \int_{\Sigma(\omega_{123})\times I^3} |A_{456}| \chi_{A_{456}}(x) ds_{456} d\omega_{456}.
\end{align*}
Given $x \in B_R^c \doteq \{y: |y|>R\}$, we define 
\[
 \Omega_\delta (x) = \{\omega\in \mathbb{S}^2: \exists s\in I, x\in A_{s,\omega,\delta}\} = \{\omega\in \mathbb{S}^2: a-\delta < x\cdot \omega<b\}.
\]
The least upper bound of $J(s_{123}, \omega_{123})$ satisfies 
\begin{align*}
 \sup_{s_{123} \in I^3, \omega_{123} \in (\mathbb{S}^2)^3} J(s_{123}, \omega_{123}) & \leq  \sup_{s_{123} \in I^3, \omega_{123} \in (\mathbb{S}^2)^3}  \sup_{x\in A_{123}} \int_{\Sigma(\omega_{123})\times I^3} |A_{456}| \chi_{A_{456}}(x) ds_{456} d\omega_{456}\\
 & \leq \sup_{x\in B_R^c, \omega_{123}\in \Omega_\delta^3(x)} \int_{(\Sigma(\omega_{123})\cap \Omega_\delta^3(x))\times I^3} |A_{456}| \chi_{A_{456}}(x) ds_{456} d\omega_{456}.
\end{align*}
We observe
\[
 |A_{456}| \leq \frac{\delta^3}{|(\omega_4\times\omega_5)\cdot \omega_6|},
\]
and obtain
\[
 \sup_{s_{123}, \omega_{123} \in I^3 \times (\mathbb{S}^2)^3} J(s_{123}, \omega_{123}) \leq \sup_{x\in B_R^c, \omega_{123}\in \Omega_\delta^3(x)} \int_{(\Sigma(\omega_{123})\cap \Omega_\delta^3(x))\times I^3} \frac{\delta^3 \chi_{A_{456}}(x)}{|(\omega_4\times\omega_5)\cdot \omega_6|} ds_{456} d\omega_{456}.
\]
Next we recall
\[
 \chi_{A_{456}} (x) = 1 \quad \Leftrightarrow \quad s_k - \delta < x\cdot \omega_k < s_k, \; \forall k = 4,5,6.
\]
Thus we have 
\begin{align*}
 \sup_{s_{123}, \omega_{123} \in I^3 \times (\mathbb{S}^2)^3} J(s_{123}, \omega_{123}) &\leq \sup_{x\in B_R^c, \omega_{123}\in \Omega_\delta^3(x)} \int_{\Sigma(\omega_{123})\cap \Omega_\delta^3(x)} \frac{\delta^6}{|(\omega_4\times\omega_5)\cdot \omega_6|} d\omega_{456}\\
 &\leq \delta^6 C_{R,I,\delta_0}
\end{align*}
for all $\delta \in (0,\delta_0)$. Here $C_{R,I,\delta_0}$ is a constant independent of $\delta \in (0,\delta_0)$
\begin{equation}
 C_{R,I,\delta_0} = \sup_{x\in B_R^c, \omega_{123}\in \Omega_{\delta_0}^3(x)} \int_{\Sigma(\omega_{123})\cap \Omega_{\delta_0}^3(x)} \frac{1}{|(\omega_4\times\omega_5)\cdot \omega_6|} d\omega_{456}
\end{equation}
Plugging this upper bound in \eqref{upper bound J delta}, we obtain 
\begin{align*}
 J_\delta  &\leq 10 C_{R,I, \delta_0}  \int_{(\mathbb{S}^2)^3 \times I^3} |G(s_1,\omega_1)|^2  |G(s_2,\omega_2)|^2  |G(s_3,\omega_3)|^2 ds_{123} d\omega_{123} \\
  &\leq 10C_{R,I,\delta_0} \|G\|_{L^2(\Rm\times \mathbb{S}^2)}^6. 
\end{align*}
We make $\delta \rightarrow 0^+$ and conclude that the following inequality holds for any small constant $\delta_0>0$.
\[
 \left\|\mathbf{T} G\right\|_{L^6(\{x\in \Rm^3: |x|>R\})}^6 \leq 10 C_{R,I,\delta_0} \|G\|_{L^2(\Rm\times \mathbb{S}^2)}^6.
\]
The remaining work is to find an upper bound of $C_{R,I,\delta_0}$. Let us first fix an $x\in B_R^c$ and determine the upper bound of 
\begin{equation} \label{def of ctx}
 C_{I, \delta_0}(x) = \sup_{\omega_{123}\in \Omega_{\delta_0}^3(x)} \int_{\Sigma(\omega_{123})\cap \Omega_{\delta_0}^3(x)} \frac{1}{|(\omega_4\times\omega_5)\cdot \omega_6|} d\omega_{456}
\end{equation}
 Without loss of generality we assume $x = (0,0,h) \in \Rm^3$. Then 
\begin{align*}
 \Omega_{\delta_0} (x) & = \left\{\omega = (x_1,x_2,x_3) \in \mathbb{S}^2: \frac{a-\delta_0}{h}< x_3< \frac{b}{h}\right\}.&  & h\geq b;\\
 \Omega_{\delta_0} (x) & = \left\{\omega = (x_1,x_2,x_3) \in \mathbb{S}^2: \frac{a-\delta_0}{h}< x_3 \leq 1\right\}, & & h\in (a-\delta_0,b);\\
 \Omega_{\delta_0} (x) & = \varnothing. & & h \leq a-\delta_0.
\end{align*}
We next apply a geometric transformation so that we may work in Euclidean space $\Rm^2$ for convenience. Let $O$ be the origin in $\Rm^3$. We consider the central projection (with center $O$) from the upper half of the sphere
\[
 \mathbb{S}_+^2 = \{(x_1,x_2,x_3): x_1^2 + x_2^2 +x_3^2 = 1, x_3>0\}
\]
to the plane $x_3=1$: (Please see figure \ref{figure central projection})
\[
 \mathbf{P}: \mathbb{S}_+^2 \rightarrow \Rm^2, \qquad Y= \mathbf{P} (x_1, x_2, x_3) = (x_1/x_3,x_2/x_3). 
\]
We have
\begin{align*}
 \Omega_{\delta_0,h}^\ast \doteq \mathbf{P} (\Omega_{\delta_0}(x) ) & = \left\{Y\in \Rm^2: \frac{\sqrt{h^2-b^2}}{b} < |Y| < \frac{\sqrt{h^2-(a-\delta_0)^2}}{a-\delta_0}\right\}.& & h\geq b;\\
 \Omega_{\delta_0,h}^\ast \doteq \mathbf{P} (\Omega_{\delta_0}(x) ) & = \left\{Y\in \Rm^2:  |Y| < \frac{\sqrt{h^2-(a-\delta_0)^2}}{a-\delta_0}\right\}, & & h\in [a-\delta_0,b);\\
 \Omega_{\delta_0,h}^\ast \doteq \mathbf{P} (\Omega_{\delta_0}(x) ) & = \varnothing, & & h< a-\delta_0;
\end{align*}
is an annulus (or a disk) and 
\[
 dY = x_3^{-3} dS(\omega).
\]
We define $Y_k = \mathbf{P} \omega_k \in \Rm^2$ and use notation $\vec{y}_k$ for the vector $OY_k \in \Rm^3$. If $\omega_k= (\omega_{k,1},\omega_{k,2},\omega_{k,3}) \in \Omega_{\delta_0} (x)$, then 
\begin{equation} \label{vector quotient}
 \vec{y}_k = \omega_{k,3}^{-1} \omega_k, \qquad \frac{h}{b} < \omega_{k,3}^{-1} < \frac{h}{a-\delta_0}.
\end{equation}
Since the distance of $O$ to the plane $x_3 =1$ is $1$, the volume $V$ of tetrahedron $OY_i Y_j Y_k$ is one third of the area of triangle $Y_i Y_j Y_k$. Thus (please see figure \ref{figure projection2}) 
\begin{align*}
 \frac{1}{3}|\triangle Y_i Y_j Y_k| = V = \frac{1}{6} \left|(\vec{y}_i \times \vec{y}_j)\cdot \vec{y}_k\right|.
\end{align*}
We may combine this with \eqref{vector quotient} and obtain
\begin{equation} \label{area ratio}
 \frac{|\triangle Y_i Y_j Y_k|}{|(\omega_i \times \omega_j)\cdot \omega_k|} = \frac{\left|(\vec{y}_i \times \vec{y}_j)\cdot \vec{y}_k\right|}{2|(\omega_i \times \omega_j)\cdot \omega_k|} = \frac{1}{2\omega_{i,3} \cdot \omega_{j,3}\cdot \omega_{k,3}} \in \left(\frac{h^3}{2b^3}, \frac{h^3}{2(a-\delta_0)^3}\right) 
\end{equation} 
Therefore we may use the reciprocal assumption on triples $\omega_{123}$ and $\omega_{456}$, as well as the assumption $b/a \leq 2$ to deduce (as long as $\delta_0>0$ is sufficiently small)
\begin{align*}
 |\triangle Y_1 Y_2 Y_3|\cdot |\triangle Y_4 Y_5 Y_6| & \geq \frac{(a-\delta_0)^6}{b^6} \max_{j_1, j_2, \cdots, j_6}  |\triangle Y_{j_1} Y_{j_2} Y_{j_3}|\cdot |\triangle Y_{j_4} Y_{j_5} Y_{j_6}|\\
 & \geq  \frac{1}{65} \max_{j_1, j_2, \cdots, j_6}  |\triangle Y_{j_1} Y_{j_2} Y_{j_3}|\cdot |\triangle Y_{j_4} Y_{j_5} Y_{j_6}|. 
\end{align*}
Here again the maximum is taken for all possible permutations of $1,2,\cdots ,6$. We still call these two triangles $\triangle Y_1 Y_2 Y_3$ and $Y_4 Y_5 Y_6$ (weakly) reciprocal to each other and use the notation 
\[
 \Sigma(Y_{123}) = \{(Y_4,Y_5,Y_6) \in (\Rm^2)^3: (Y_1,Y_2,Y_3), (Y_4,Y_5,Y_6)\; \hbox{are reciprocal}\}.
\]
We apply change of variables on the integral in \eqref{def of ctx}, utilize \eqref{area ratio} and obtain
\begin{align*}
 C_{I,\delta_0}(x) & = \sup_{\omega_{123}\in \Omega_{\delta_0}^3(x)} \int_{\mathbf{P}(\Sigma(\omega_{123})\cap \Omega_{\delta_0}^3(x))}\frac{\omega_{4,3}^2 \omega_{5,3}^2 \omega_{6,3}^2}{2}\cdot \frac{1}{|\triangle Y_4 Y_5 Y_6|} dY_{456}\\
 & \leq  \frac{b^6}{2h^6} \sup_{Y_{123} \in (\Omega_{\delta_0,h}^\ast)^3}   \int_{\Sigma(Y_{123})\cap (\Omega_{\delta_0,h}^\ast)^3} \frac{1}{|\triangle Y_4 Y_5 Y_6|} dY_{456}.
\end{align*}
In summary we have 

\begin{lemma} \label{transform to geometric problem}
Assume that $b>a>0$ with $b/a\leq 2$. Let $G(s,\omega) \in L^2(\Rm\times \mathbb{S}^2)$ be supported in the region $[a,b]\times \mathbb{S}^2$. Then the function 
\[
 \mathbf{T} G (x) = \int_{\mathbb{S}^2} G(x\cdot \omega, \omega) d\omega, \qquad x\in \Rm^3
\]
satisfies the following inequality for all sufficiently small $\delta>0$: 
\[
 \int_{|x|>R} |\mathbf{T} G (x)|^6 dx \leq \left(\sup_{h>\max\{R,a-\delta\}} C_{a,b,\delta}(h) \right)\|G\|_{L^2(\Rm\times \mathbb{S}^2)}^6.
\]
The constant $C_{a,b,\delta}(h)$ is defined by 
\[
 C_{a,b,\delta}(h) = \frac{5b^6}{h^6} \sup_{Y_1,Y_2,Y_3 \in \Omega_{\delta,h}^\ast} \int_{\Sigma(Y_1Y_2Y_3)\cap (\Omega_{\delta,h}^\ast)^3} \frac{1}{|\triangle Y_4 Y_5 Y_6|} dY_4 dY_5 dY_6.
\]
Here $\Omega_{\delta,h}^\ast$ is an annulus (or disk) region in $\Rm^2$ defined by
\begin{align*}
 \Omega_{\delta,h}^\ast &\doteq  \left\{Y\in \Rm^2: \frac{\sqrt{h^2-b^2}}{b} < |Y| < \frac{\sqrt{h^2-(a-\delta)^2}}{a-\delta}\right\}, & & h\geq b;\\
 \Omega_{\delta,h}^\ast &\doteq  \left\{Y\in \Rm^2:  |Y| < \frac{\sqrt{h^2-(a-\delta)^2}}{a-\delta}\right\}, & & h\in (a-\delta, b).
\end{align*}
And $\Sigma(Y_1 Y_2 Y_3)$ consists of all (weakly) reciprocal triples of $(Y_1, Y_2,Y_3)$ in $\Rm^2$:
\begin{align*}
 \Sigma (Y_1 Y_2 Y_3) = \left\{(Y_4, Y_5, Y_6): |\triangle Y_1 Y_2 Y_3|\cdot |\triangle Y_4 Y_5 Y_6| \geq \frac{1}{65} \max_{j_1, j_2, \cdots, j_6}  |\triangle Y_{j_1} Y_{j_2} Y_{j_3}|\cdot |\triangle Y_{j_4} Y_{j_5} Y_{j_6}|\right\}.
\end{align*}
Here the maximum is taken for all possible permutations of $1,2,\cdots, 6$. 
\end{lemma}

\section{Geometric Observations}

In this section we make some geometric observations. We first give a few geometric characteristics of (weakly) reciprocal triangles in $\Rm^2$ and then a few properties an annulus region satisfies. Many of the following results are simple geometric observations and might have been previously known. Here we still give their proof for the reason of completeness. In this section we say that a triangle $\triangle ABC$ is of size $L$ if and only if $L \leq \max\{|AB|, |BC|, |CA|\} < 2L$. 

\subsection{Reciprocal triangles}
\begin{figure}[h]
 \centering
 \includegraphics[scale=1.25]{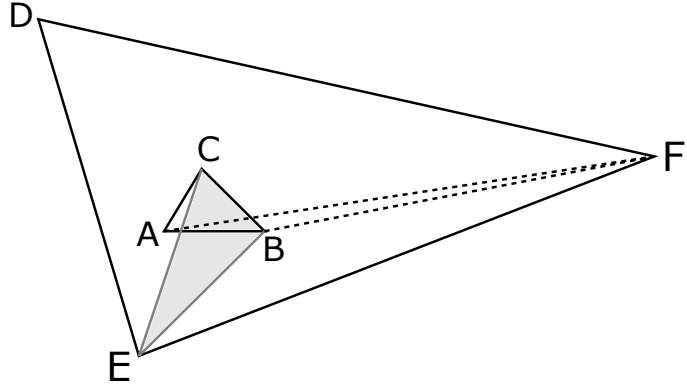}
 \caption{Illustration of case 2} \label{figure reciprocal_case2}
\end{figure}

In this subsection, we consider (weakly) reciprocal triangles in $\Rm^2$, as defined in the previous section. 
\begin{lemma} \label{parallelogram law}
 Let $\triangle ABC$ be of size $L$ and $D\in \Rm^2$ satisfy $d = d(D, \triangle ABC)\gg L$. Then either $|\triangle DAB| \gtrsim (d/L) |\triangle ABC|$ or $|\triangle DAC| \gtrsim (d/L) |\triangle ABC|$. 
\end{lemma}
\begin{proof}
We always have $\max\{\sin \angle DAC, \sin \angle DAB\} \geq (1/2) \sin \angle BAC$. Thus
\begin{align*}
 \max\{|\triangle DAC|, |\triangle DAB|\} &\gtrsim \max\{|DA| \cdot |AC| \sin \angle DAC, |DA|\cdot |AB| \sin \angle DAB\}\\
 & \gtrsim (d/L) |AB|\cdot |AC| \max\{\sin \angle DAC, \sin \angle DAB\} \\
 & \gtrsim (d/L) |AB|\cdot |AC| \sin \angle BAC \\
 & \gtrsim (d/L) |\triangle BAC|.
\end{align*} 
\end{proof}
\noindent This immediately gives
\begin{corollary} \label{parallelogram law 2}
 Let $\triangle ABC$ be of size $L$ and $D\in \Rm^2$ satisfy $d = d(D, \triangle ABC)\gg L$. Then at least two of the following inequalities holds
\begin{align*}
 &|\triangle DAB| \gtrsim (d/L) |\triangle ABC|;& &|\triangle DBC| \gtrsim (d/L) |\triangle ABC|;& &|\triangle DCA| \gtrsim (d/L) |\triangle ABC|.&
\end{align*}
\end{corollary} 
\begin{proposition} \label{close to vertices}
 Let $\triangle ABC, \triangle DEF \subset \Rm^2$ be reciprocal and of sizes $L \ll M$, respectively. Then there exists a vertex of $\triangle DEF$ (say $D$) so that $|AD|, |BD|, |CD| \lesssim L$.  
\end{proposition}
\begin{proof}
Let us prove Proposition \ref{close to vertices} by contradiction.  We assume 
 \[
  |AD|, |AE|, |AF|, |BD|, |BE|, |BF|, |CD|, |CE|, |CF| \gg L.
 \]
Without loss of generality we also assume $|DF| \geq |EF| \geq |DE|$. Thus $|DF|, |EF|\simeq M$. We consider two cases: case 1, $\triangle ABC$ is close to the vertex $F$; case 2, $\triangle ABC$ is far away from the vertex $F$. 
 
 \paragraph{Case 1}  If $|AF|, |BF|, |CF| \ll M$. We apply Corollary \ref{parallelogram law 2} on $\triangle ABC$ and $F$, at least two of the following holds 
\begin{align*}
 &|\triangle FAB| \gg |\triangle ABC|;& &|\triangle FBC| \gg  |\triangle ABC|;& &|\triangle FCA| \gg |\triangle ABC|.&
\end{align*}
Similarly at least two of the following inequalities holds
\begin{align*}
 &|\triangle EAB| \gg |\triangle ABC|;& &|\triangle EBC| \gg  |\triangle ABC|;& &|\triangle ECA| \gg |\triangle ABC|.&
\end{align*}
Thus we may find two vertices from $\triangle ABC$, say $AB$, so that we have 
\[
 |\triangle FAB|, |\triangle EAB| \gg |\triangle ABC|.
\]
Next we show that either $|\triangle CDE| \gtrsim |\triangle DEF|$ or $|\triangle CDF| \gtrsim |\triangle DEF|$ holds. This immediately gives a contradiction to our reciprocal assumption. In fact, if the first inequality fails, i.e. $|\triangle CDE| \ll |\triangle DEF|$, then we have 
\[
 d(C, DE) \ll d(F, DE). 
\]
Our assumption $|CF| \ll M$ guarantees that $|CD| \simeq M \simeq |DF|$, thus we have 
\[
 \sin \angle CDE = \frac{d(C,DE)}{|CD|} \ll \frac{d(F,DE)}{|DF|} = \sin \angle FDE.
\]
It immediately follows that $\sin \angle FDC \simeq \sin \angle FDE$. Thus 
\[
 |\triangle CDF| = |CD| \cdot |DF| \sin \angle FDC \gtrsim |DE| \cdot |DF| \sin \angle FDE = |\triangle DEF|. 
\]
This finishes the argument in case one. Please see figure \ref{figure reciprocal_case1} for an illustration of the proof. 

 \begin{figure}[h]
 \centering
 \includegraphics[scale=1.25]{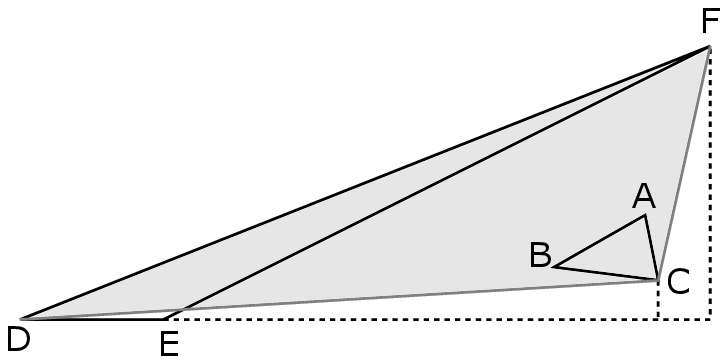}
 \caption{Illustration of case 1} \label{figure reciprocal_case1}
\end{figure}

\paragraph{Case 2} In this case $|AF|, |BF|, |CF| \gtrsim M$. Given any vertex $X\in \{A,B,C\}$, we have either $\sin \angle XFE \geq (1/2) \sin \angle DFE$ or $\sin \angle XFD \geq (1/2) \sin \angle DFE$. As a result, we may find one vertex from $\{D,E\}$ (say $D$) and two vertices from $\{A,B,C\}$ (say $A,B$) so that 
\begin{align*}
 &\sin \angle AFD \geq (1/2) \sin \angle DFE;& &\sin \angle BFD \geq (1/2) \sin \angle DFE.&
\end{align*}
Combining these angles with our assumptions $|AF|, |BF|\gtrsim M$ and $|DF|, |EF| \simeq M$, we obtain
\begin{align} \label{two points bigger}
 &|\triangle AFD| \gtrsim |\triangle DEF|;& &|\triangle BFD| \gtrsim |\triangle DEF|.& 
\end{align}
Finally we apply Lemma \ref{parallelogram law} on $\triangle CAB$ and $E$ to conclude that either $|\triangle EBC| \gg |\triangle ABC|$ or $|\triangle ECA| \gg |\triangle ABC|$ holds. A combination of this with \eqref{two points bigger} immediately gives a contradiction. Please see figure \ref{figure reciprocal_case2} for an illustration of this case. Combining case 1 and 2, we finish the proof of Proposition \ref{close to vertices}. 
\end{proof} 

\begin{corollary} \label{close to each other}
 Let $\triangle ABC, \triangle DEF \subset \Rm^2$ be reciprocal of sizes $L, M$, respectively. Then they can not be too far away from each other. Namely we always have
 \[
  d(\triangle ABC, \triangle DEF) \lesssim \min\{L, M\}. 
 \]
\end{corollary}
\begin{proof}
This corollary clearly holds if the size of one triangle is much larger than that of the other, thanks to Proposition \ref{close to vertices}. Thus we only need to consider the case $L \simeq M$. If the corollary failed, we would have 
\[
 d(\triangle ABC, \triangle DEF) \gg L, M. 
\]
We may apply Corollary \ref{parallelogram law 2} on the triangle $DEF$ and the point $A$, then on the same triangle and the point $B$. This enable us to find two vertices from $DEF$ (say $DE$) so that 
\begin{align*}
 &|\triangle ADE| \gg |\triangle DEF|;& &|\triangle BDE| \gg |\triangle DEF|.&
\end{align*}
We then apply Lemma \ref{parallelogram law} on the triangle $CAB$ and the point $F$, then conclude that at least one of the following holds
\begin{align*}
 &|\triangle BCF| \gg |\triangle ABC|;& &|\triangle ACF| \gg |\triangle ABC|.&
\end{align*}
Either of these contradicts our reciprocal assumption. 
\end{proof}

\begin{figure}[h]
 \centering
 \includegraphics[scale=1.25]{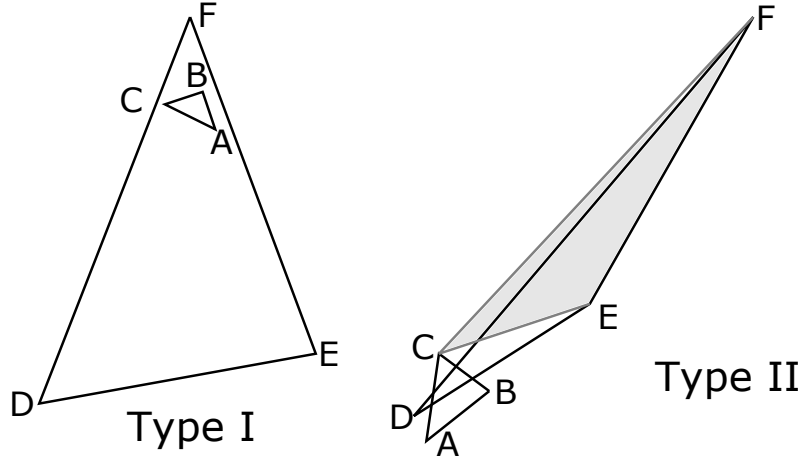}
 \caption{Classification of reciprocal triangles} \label{figure reciprocal_class}
\end{figure}

\begin{proposition}[Classification]
Let $\triangle ABC$ and $\triangle DEF$ be two reciprocal triangles of sizes $L\ll M$, respectively. Without loss of generality we also assume 
$|BC|$ and $|DE|$ are the shortest edge in the corresponding triangles. Then the location of smaller triangle $ABC$ satisfies either of the following 
\begin{itemize}
 \item[(I)] $|AF|, |BF|, |CF| \lesssim L$;
 \item[(IIa)] $|AD|, |BD|, |CD| \lesssim L$ so that $\max\{|\triangle BEF|, |\triangle CEF|\} \gtrsim |\triangle DEF|$; 
 \item[(IIb)] $|AE|, |BE|, |CE| \lesssim L$ so that $\max\{|\triangle BDF|, |\triangle CDF|\} \gtrsim |\triangle DEF|$.
\end{itemize}
We call these triangles Type I reciprocal if they satisfies (I) and call them Type II reciprocal if they satisfies either (IIa) or (IIb). Please see figure \ref{figure reciprocal_class}.
\end{proposition}
\begin{proof}
Proposition \ref{close to vertices} guarantees that if (I) fails, then we have either $|AD|, |BD|, |CD|\lesssim L$ or $|AE|, |BE|, |CE| \lesssim L$. Without loss of generality, we assume $|AD|, |BD|, |CD| \lesssim L$ and show that either (IIa) or (IIb) holds. Because $|FB|, |FD|, |FE| \simeq M$, we may conclude that either $|\triangle BFD| \gtrsim |\triangle DEF|$ or $|\triangle BEF| \gtrsim |\triangle DEF|$ holds by considering the angles $\angle BFD$ and $\angle BFE$. If the latter holds, $\triangle ABC$ satisfies (IIa). Thus we only need to consider the first case. Similarly we may assume $|\triangle CFD| \gtrsim |\triangle DEF|$. Now we claim that $|AE|, |BE|, |CE| \lesssim L$ thus (IIb) holds. Otherwise we may apply Lemma \ref{parallelogram law} and conclude that either $|\triangle ABE| \gg |\triangle ABC|$ or $|\triangle ACE| \gg |\triangle ABC|$. This means
\[
 \max\{|\triangle BFD| \cdot |\triangle ACE|, |\triangle CFD| \cdot |\triangle ABE|\} \gg |\triangle ABC| \cdot |\triangle DEF|,
\]
thus contradicts the reciprocal assumption. 
\end{proof}

\subsection{About annulus}

In this subsection we give a few geometric properties of a circular annulus region. We consider a circular annulus region $\Omega^\ast \subset \Rm^2$, whose outer radius is $r$, inner radius is $r-w$ and width is $w$. We will also use the notation $O$ for the center of $\Omega^\ast$. 

\begin{lemma} \label{range of angle 1}
 Assume $A,B\in \Omega^\ast$ and $AC \perp OA$. Then 
\[
 \sin \angle BAC \leq \max\left\{\frac{2w}{|AB|}, \frac{2|AB|}{r}\right\}.
\]
\end{lemma}
 \begin{figure}[h]
 \centering
 \includegraphics[scale=1.25]{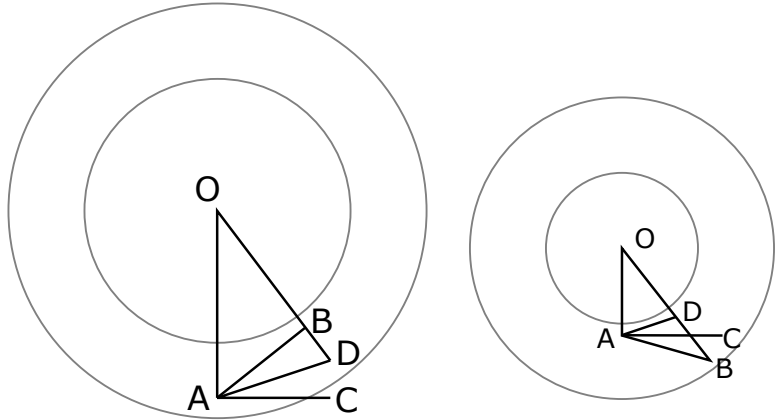}
 \caption{Illustration of proof} \label{figure annulus1}
\end{figure}
\begin{proof}
 First of all, if $|AB|\geq r/2$ or $|AB| \leq 2w$, then the right hand side is greater or equal to $1$, thus the inequality holds. We now assume $2w < |AB| < r/2$ thus $w<r/4$. Let $D$ be the point on the ray $OB$ so that $|OD| = |OA|$. We have 
  \[
  \sin \angle CAD = \frac{|AD|}{2|OA|} \leq \frac{|AB|+w}{2(r-w)} \leq \frac{3|AB|/2}{3r/2} = \frac{|AB|}{r}. 
 \] 
  We also have
 \[
  \sin \angle BAD = \frac{|BD| \sin \angle BDA}{|AB|} \leq \frac{w}{|AB|}.
 \]
 Finally we have
 \[
  \sin \angle BAC \leq \sin \angle CAD + \sin \angle BAD \leq \frac{|AB|}{r} + \frac{w}{|AB|} \leq \max\left\{\frac{2w}{|AB|}, \frac{2|AB|}{r}\right\}.
 \]
\end{proof}

\begin{corollary} \label{range of angle 2}
 Assume $L \leq r$ and $A \in \Omega^\ast$. Then 
 \begin{itemize}
  \item[(a)] $\left|\left\{\Theta \in \mathbb{S}^1: \exists l \in [L/2,2L], \, \hbox{s.t.}\, A+ l\Theta \in \Omega^\ast\right\}\right| \leq 8\pi \max\{w/L, L/r\}$. Here $A+l \Theta$ is the terminal point of the vector in $\Rm^2$ with starting point $A$, length $l$ and direction $\Theta$. 
  \item[(b)] If $B,C \in \Omega^\ast$ so that $L/2\leq |AB|, |AC| \leq 2L$, then we have 
\[
 \sin \angle BAC \leq 8 \max\left\{w/L, L/r\right\}
\]
 \end{itemize}
\end{corollary}
\begin{proof}
 Let $AD \perp OA$ and $E=A+l \Theta \in \Omega^\ast$, $l \in [L/2,2L]$. By Lemma \ref{range of angle 1}, we have 
 \[
  \sin \angle EAD \leq 4 \max\left\{w/L, L/r\right\}.
 \]
 We observe ($z\in [0,1]$)
 \begin{align*}
  &\sin \angle EAD \leq z \quad \Leftrightarrow \angle EAD \in [0,\arcsin z]\cup [\pi -\arcsin z, \pi];& &\arcsin z \leq \pi z/2.&
 \end{align*}
 Thus the subset of $\mathbb{S}^1$ consisting all possible directions of $AE$ has a measure smaller or equal to $8\pi \max\{w/L, L/r\}$. This proves part (a). For part (b), a similar argument shows
 \[
  \sin \angle DAB, \sin \angle DAC \leq 4 \max\left\{w/L, L/r\right\}.
 \]
 Thus 
 \[
  \sin \angle BAC \leq \sin \angle DAB + \sin \angle DAC \leq 8 \max\left\{w/L, L/r\right\}.
 \]
\end{proof}

\begin{figure}[h]
 \centering
 \includegraphics[scale=1.25]{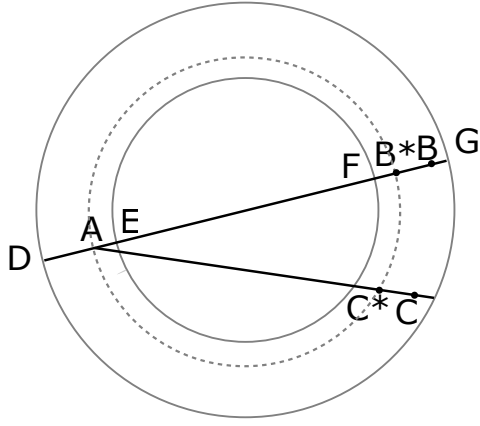}
 \caption{estimate of $|BC|$} \label{figure annulus4}
\end{figure}

\begin{lemma} \label{length of edge}
 Let $A, B, C\in \Omega^\ast$ so that $|AB|, |AC|\geq 3\sqrt{wr}$. Then we have 
\[
 2r \sin \angle BAC - 2\sqrt{wr} -2w < |BC| < 2r \sin \angle BAC + 2\sqrt{wr}.
\]
\end{lemma}
\begin{proof}
First of all, we claim that the line $AB$ must intersect the inner boundary of $\Omega^\ast$ at two different points, otherwise the length $|AB|$ can never exceed $2\sqrt{w(2r-w)}$. Let $D,E,F,G$ be the intersection points of the line $AB$ with the boundary of $\Omega^\ast$, as shown in figure \ref{figure annulus4}, so that $A$ is on the line segment $DE$. We have $|DG| > |AB| \geq 3\sqrt{wr}$. In addition
 \[
  (|DG|-|FG|)\cdot |FG| = |DF|\cdot |FG| = w(2r-w) < 2wr, \qquad |FG|\leq |DG|/2.
 \]
 This immediately gives $|DE| = |FG| < \sqrt{wr}$. As a result, $B$ must be on the line segment $FG$. Let $B^\ast$ be the point on the line segment $FG$ so that $|OA| = |OB^\ast|$. We always have $|BB^\ast| < |FG| < \sqrt{wr}$. We may define $C^\ast$ in a similar way, as shown in figure \ref{figure annulus4}. Again we have $|CC^\ast| < \sqrt{wr}$. Since $A, B^\ast, C^\ast$ is on the same circle of radius $|OA| \in (r-w,r)$, we have 
\[
 2(r-w) \sin \angle BAC < |B^\ast C^\ast| < 2r \sin \angle BAC.
\]
Therefore 
\begin{align*}
|BC| & \leq |BB^\ast| + |B^\ast C^\ast| + |C^\ast C| < 2r \sin \angle BAC + 2\sqrt{wr};\\
|BC| & \geq |B^\ast C^\ast| - |B B^\ast| -  |C^\ast C| > 2(r-w) \sin \angle BAC - 2\sqrt{wr} \geq 2r\sin \angle BAC - 2\sqrt{wr} -2w.
\end{align*}
\end{proof}
\begin{corollary} \label{big triangle in omega}
 Let $A, B, C \in \Omega^\ast$ so that $|AB|, |BC|, |CA| \geq 4\sqrt{wr}$. Then
 \begin{itemize}
  \item[(a)] $r \sin \angle BAC < |BC| < 4r \sin \angle BAC$;
  \item[(b)] $|\triangle ABC| \simeq |AB|\cdot |BC| \cdot |CA|/r$. 
 \end{itemize}
\end{corollary}
\begin{proof}
 We may rewrite the conclusion of Lemma \ref{length of edge} in the form of
 \[
  |BC| - 2\sqrt{wr} < 2r \sin \angle BAC < |BC| + 2\sqrt{wr} + 2w.
 \]
 We then combine this inequality with the assumption $|BC| \geq 4\sqrt{wr}$
 \[
  \frac{1}{2} |BC| < 2r \sin \angle BAC < 2|BC|.
 \]
This proves part (a). Part (b) immediately follows part (a) and the basic formula 
 \[
  |\triangle ABC| = \frac{1}{2} |AB| \cdot |AC|\sin \angle BAC.
 \]
 
\end{proof}
\begin{corollary} \label{semi-big triangle in omega}
 Let $A, B, C \in \Omega^\ast$ so that $|AB|, |AC| \geq 4\sqrt{wr}$. Then 
\[
 |\triangle ABC| \lesssim \frac{|AB|\cdot |AC|\cdot \max\{|BC|, \sqrt{wr}\}}{r}.
\]
\end{corollary}
\begin{proof}
 If $|BC| \geq 4\sqrt{wr}$, then we may apply Corollary \ref{big triangle in omega} and finish the proof. If $|BC| < 4\sqrt{wr}$, then Lemma \ref{length of edge} implies 
 \[
  2r \sin BAC < |BC| + 2\sqrt{wr} + 2w < 8 \sqrt{wr} \quad \Rightarrow \quad \sin BAC < 4 \sqrt{wr} /r.
 \]
 This immediately gives
 \[
  |\triangle ABC| =  \frac{1}{2} |AB| \cdot |AC|\sin \angle BAC \leq \frac{|AB|\cdot |AC|\cdot 2\sqrt{wr}}{r}.
 \]
\end{proof}

\begin{figure}[h]
 \centering
 \includegraphics[scale=1.25]{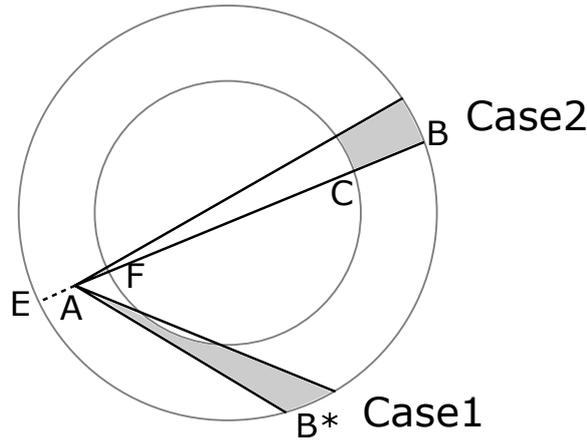}
 \caption{Area by angle} \label{figure annulus2}
\end{figure}

\begin{lemma}[Area by angle] \label{area upper bound by angle}
 Let $A \in \Omega^\ast$ and $K \subset \mathbb{S}^1$ be measurable. Then 
 \[
 \left|\Omega^\ast \cap \{A+l \Theta\in \Rm^2: l\in \Rm^+, \Theta \in K\}\right| \leq 4wr |K|.
 \]
\end{lemma}
\begin{proof} 
 It suffices to consider the case $K = (\theta, \theta+d \theta)$. Here we slightly abuse the notation, the angle $\theta$ actually represent the direction $\Theta = (\cos\theta, \sin \theta) \in \mathbb{S}^1$. Let $B$ (or $B^\ast$) be the point where the ray $A+l \Theta$ meets the outer boundary of the annulus $\Omega^\ast$. We consider two cases. Case 1, if $|AB^\ast| \leq 2\sqrt{2wr}$ is relatively short, then we have 
 \[
  dS \leq \frac{1}{2}|AB^\ast|^2 d\theta \leq 4 wr d\theta. 
 \]
 Case 2, if $|AB| > 2 \sqrt{2wr}$ is long, then we claim that the segment $AB$ must intersect with the inner boundary of $\Omega^\ast$ at two different points. Otherwise the length $|AB|$ can never exceed $2\sqrt{w(2r-w)}$. Let $E, F, C, B$ be the intersection points of line $AB$ with the boundary circles of $\Omega^\ast$, as shown in figure \ref{figure annulus2}. We have 
 \[
  |EC|\cdot |BC| = w(2r-w).
 \]
 Thus we have ($|AF|\leq |EF| = |BC| \leq |EC|$)
 \begin{align*}
  dS = \left[\left(|AC|+\frac{1}{2}|BC|\right)|BC|+\frac{1}{2}|AF|^2\right] d\theta \leq \left(|EC| + |EF|\right) |BC| d\theta \leq 4 wr d\theta.
 \end{align*}
\end{proof}

\begin{corollary} \label{area upper bound by sine}
  Let $A \in \Omega^\ast$, $B \in \Rm^2\setminus \{A\}$ and $z \in \Rm^+$. Then 
 \[
  \left|\{C\in \Omega^\ast: \sin \angle BAC \leq  z \}\right| \leq 8\pi zwr. 
 \]
\end{corollary}
\begin{proof}
 If $z\geq 1$, then the inequality is trivial since $|\Omega^\ast| \leq 2\pi wr$. If $z\in (0,1)$, then 
\[
 \sin \angle BAC \leq z \quad \Leftrightarrow \quad \angle BAC \in [0, \arcsin z] \cup [\pi-\arcsin z, \pi]. 
\]
We then utilize the inequality $\arcsin z \leq \pi z/2$ and apply Lemma \ref{area upper bound by angle} to complete the proof. 
\end{proof}

\begin{remark} \label{trivial area}
The following will also be used in the subsequent section: Assume $A, B\in \Rm^2$, $L, z\in \Rm^+$. Let $K \subset \mathbb{S}^1$ be measurable. Then we have
 \begin{align*} 
  \left|\{A+l \Theta\in \Rm^2: l\in (0,L), \Theta \in K\}\right| & \leq \frac{1}{2} L^2 |K|.  \\
  \left|\{C\in \Rm^2: |CA| \leq L,  \sin \angle BAC \leq  z \}\right| & \leq \pi L^2 z.
 \end{align*}
\end{remark}

\begin{figure}[h]
 \centering
 \includegraphics[scale=1.25]{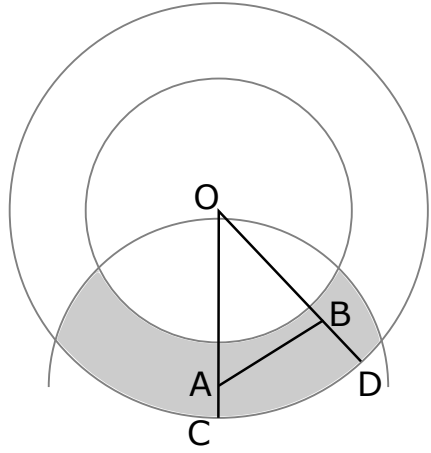}
 \caption{Area by distance} \label{figure annulus3}
\end{figure}

\begin{lemma}[Area by distance] \label{area upper bound by distance}
 Let $A \in \Omega^\ast$ and $L>0$. Then 
\[
 \left|\Omega^\ast \cap B(A,L)\right|  \leq 2\pi Lw.
\]
Here $B(A,L)$ is the disk of radius $L$ centered at $A$. 
\end{lemma}
\begin{proof}
 This is trivial if $L<2w$ because in this case $2\pi Lw > \pi L^2 = |B(A,L)|$. Let us assume $L \geq 2w$. Given any point $B\in B(A,L) \cap \Omega^\ast$, let $C,D$ be the intersection points of the rays $OA$, $OB$ with the outer boundary of $\Omega^\ast$, as shown in figure \ref{figure annulus3}. We have
\[
 2r \sin \frac{\angle AOB}{2} = |CD| \leq |AB| + |AC| +|BD| \leq L + 2w \leq 2L.
\]
Thus $\angle AOB \leq \pi L/r$. This immediately gives 
\[
  \left|\Omega^\ast \cap B(A,L)\right| \leq \frac{2\pi L}{r} \cdot wr = 2\pi Lw.
\]
\end{proof}

\section{Proof of Geometric Inequality}
In this section we prove 

\begin{proposition} \label{geometric inequality}
 Let $\Omega^\ast \subset \Rm^2$ be a circular annulus region with outer radius $r$ and width $w$. Here $w \leq r$. Then
 \[
  \sup_{D,E,F \in \Omega^\ast} \int_{\Sigma(DEF)\cap (\Omega^\ast)^3} \frac{1}{|\triangle XYZ|} dX dY dZ \lesssim w^3 r.
 \]
 Here $\Sigma(DEF)$ is the set of all reciprocal triples of $(D,E,F)$ in $\Rm^2$, as defined in Lemma \ref{transform to geometric problem}. 
\end{proposition}

\begin{figure}[h]
 \centering
 \includegraphics[scale=1.25]{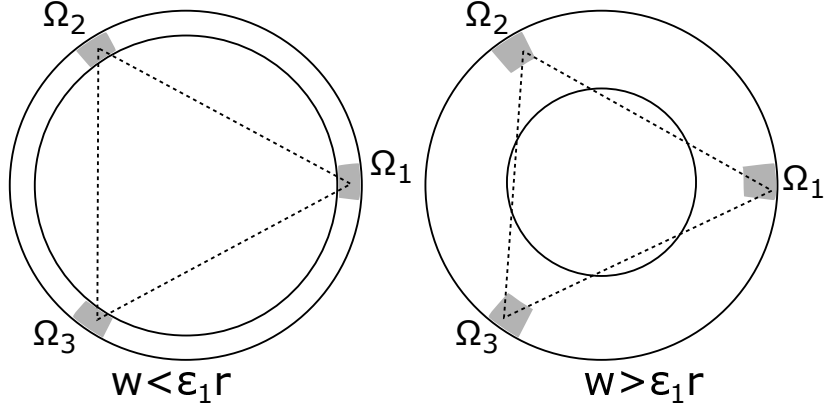}
 \caption{Optimal upper bound} \label{figure optimalbound}
\end{figure}

\begin{remark}
 The upper bound given above is optimal. We choose three angles $\theta_1 = 0$, $\theta_2 = 2\pi/3$, $\theta_3 = 4\pi/3$, and three regions accordingly by polar coordinates ($\eps_1$ is a small constant)
\[
 \Omega_k = \{(\rho\cos \theta, \rho \sin \theta): r- \min\{w, \eps_1 r\}<\rho < r, \theta_k -\eps_1 < \theta < \theta_k +\eps_1\},
\]
as show in figure \ref{figure optimalbound}. If we choose triples $(D,E,F), (X,Y,Z) \in \Omega_1\times \Omega_2\times \Omega_3$, then $\triangle DEF$ and $\triangle XYZ$ are reciprocal to each other, as long as the constant $\eps_1$ is sufficiently small. It is because these triangles are among the biggest triangles in the disk of radius $r$. This implies if we fix $(D,E,F) \in \Omega_1 \times \Omega_2 \times \Omega_3$, then 
\[
 \int_{\Sigma(DEF)\cap (\Omega^\ast)^3} \frac{dX dY dZ}{|\triangle XYZ|} \gtrsim \int_{\Omega_1 \times \Omega_2\times \Omega_3} \frac{dX dY dZ}{|\triangle XYZ|} \gtrsim \frac{(\eps_1 r \min\{w, \eps_1 r\})^3}{r^2} \gtrsim w^3 r.
\]
\end{remark}

\paragraph{Sizes and angles} In order to take advantage of the geometric properties of reciprocal triangles, we sort all reciprocal triangles $\triangle XYZ$ by their sizes and angles. We choose dyadic sequences of sizes:
\[
 L \in \{r,r/2,r/4, \cdots\}. 
\]
We say that $\triangle XYZ$ is of size $L$ if and only if $L \leq \max \{|XY|, |XZ|, |YZ|\} < 2L$. Without loss of generality we also assume that $\angle YXZ$ is the smallest among the three angles of $\triangle XYZ$. Thus we have $|YZ| = \min\{|XY|, |XZ|, |YZ|\}$. If $\triangle XYZ$ is of size $L$, then $L/2\leq |XY|, |XZ|<2L$. As a result, we define (the upper bound of $\phi_L$ can be determined by Lemma \ref{range of angle 2})
\[
 \phi_L \doteq \sup \{\sin \angle YXZ: X,Y,Z\in \Omega^\ast, L/2 \leq |XY|, |XZ|<2L\} \lesssim \max\{w/L, L/r\}. 
\]
and 
\[
 \Phi_n^L = \{\theta\in (0,\pi): 2^{-n-1} \phi_L < \sin \theta \leq 2^{-n} \phi_L\}, \qquad n\geq 0.
\]
We always have 
\begin{equation} \label{measure of phi}
 \left|\Phi_n^L\right| \lesssim 2^{-n} \phi_L.
\end{equation}
We then sort all reciprocal triangles $\triangle XYZ$ of a given triangle $\triangle DEF$ by their sizes and angles. We define 
\[
 \Omega_{L,n} = \left\{(X,Y,Z)\in (\Omega^\ast)^3: \begin{array}{c} \triangle XYZ\, \hbox{is a reciprocal triangle of} \, \triangle DEF \,
  \hbox{whose size is $L$} \\ \hbox{and whose smallest interior angle $\angle YXZ$ is in} \,\Phi_n^L\end{array} \right\}. 
\]
We immediately have for a fixed triangle $\triangle DEF$
\begin{equation}
 \int_{\Sigma(DEF)\cap (\Omega^\ast)^3} \frac{1}{|\triangle XYZ|} dX dY dZ \leq 3\sum_{L, n} \int_{\Omega_{L,n}} \frac{1}{|\triangle XYZ|} dX dY dZ. 
\end{equation}
For convenience we also assume that the size of $\triangle DEF$ is $M$ and the smallest angle of $\triangle DEF$ is $\angle EFD$. We split the big sum in the right hand side into three parts: large sizes $L\gg M$, small sizes $L \ll M$ and comparable sizes $L \simeq M$. 

\subsection{Large sizes}
We first consider the case that the size $L$ of $\triangle XYZ$ is much larger than that of $\triangle DEF$. According to our classification of reciprocal triangles, we consider two cases, i.e. Type I reciprocal triangles and Type II reciprocal triangles. We write 
\[ 
 \Omega_{L,n} = \Omega_{L,n}^1 \cup \Omega_{L,n}^2.
\]
Here 
\begin{align*}
 \Omega_{L,n}^1  & =  \left\{(X,Y,Z)\in \Omega_{L,n}: \triangle XYZ \; \hbox{and} \; \triangle DEF \; \hbox{are Type I reciprocal} \right\};\\
 \Omega_{L,n}^2  & =  \left\{(X,Y,Z)\in \Omega_{L,n}: \triangle XYZ \; \hbox{and} \; \triangle DEF \; \hbox{are Type II reciprocal} \right\}.
\end{align*}

\begin{figure}[h]
 \centering
 \includegraphics[scale=1.25]{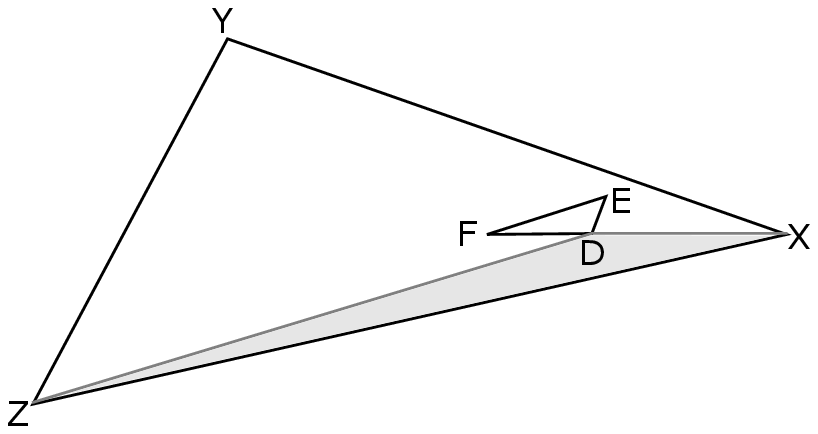}
 \caption{Large size, Type I reciprocal triangles} \label{figure large_type1}
\end{figure}

\paragraph{Type I} In this case we have $|DY|, |EY|, |FY| \simeq L \gg M$. According to Lemma \ref{parallelogram law}, we have either $|\triangle FDY| \gtrsim (L/M) |\triangle DEF|$ or $|\triangle FEY| \gtrsim (L/M) |\triangle DEF|$. Without loss of generality let us assume the latter one.\footnote{Strictly speaking, we need to consider both two cases. The argument here only takes care of one case. The other case can be dealt with in exactly the same way.}  A combination of this and the reciprocal assumption implies 
\[
 |\triangle DXZ| \lesssim (M/L) |\triangle XYZ|\quad \Rightarrow \quad |DX|\cdot |XZ| \sin \angle DXZ \lesssim (M/L) |XY|\cdot |XZ| \sin \angle YXZ.
\]
Thus we have 
\[
 |DX| \sin \angle DXZ \lesssim M \sin \angle YXZ \simeq M \cdot 2^{-n} \phi_L.
\]
This means that if $(X,Y,Z) \in \Omega_{L,n}^1$, then at least one of the following holds (see figure \ref{figure large_type1})
\begin{itemize}
 \item $|DX| \lesssim 2^{-n/2} r^{1/2} M^{1/2} \phi_L^{1/2}$; 
 \item $\sin \angle DXZ \lesssim 2^{-n/2} r^{-1/2} M^{1/2} \phi_L^{1/2}$.
\end{itemize}
We may write $\Omega_{L, n}^1 = \Omega_{L,n}^{1,1} \cup \Omega_{L,n}^{1,2}$ as a union of two parts accordingly. Here we define 
\begin{align*}
 \Omega_{L,n}^{1,1} & = \left\{(X,Y,Z)\in \Omega_{L,n}^{1}: |DX| \lesssim 2^{-n/2} r^{1/2} M^{1/2} \phi_L^{1/2} \right\};\\
 \Omega_{L,n}^{1,2} & = \left\{(X,Y,Z)\in \Omega_{L,n}^{1}: \sin \angle DXZ \lesssim 2^{-n/2} r^{-1/2} M^{1/2} \phi_L^{1/2} \right\}.
\end{align*} 
Now we are ready to find the upper bounds of the integrals ($k =1,2$)
\[
 \int_{\Omega_{L,n}^{1,k}} \frac{1}{|\triangle XYZ|} dX dY dZ. 
\]
Let us first consider the case $L\geq \sqrt{wr}$ and $k=1$. It is clear that 
\[
 \int_{\Omega_{L,n}^{1,1}} \frac{1}{|\triangle XYZ|} dX dY dZ \lesssim \frac{\left|\Omega_{L,n}^{1,1} \right|}{2^{-n} L^2 \phi_L}. 
\]
Next we give an upper bound of the measure of $\Omega_{L,n}^{1,1}$. First of all, we observe 
\[
 \Omega_{L,n}^{1,1} \subset \tilde{\Omega} \doteq \left\{(X,Y,Z)\in (\Omega^\ast)^3: |DX| \lesssim 2^{-n/2} r^{1/2} M^{1/2} \phi_L^{1/2}, |XY|<2L, \angle ZXY \in \Phi_n^L\right\}. 
\]
Thus we may find an upper bound of the measure of $\tilde{\Omega}$ instead. According to Lemma \ref{area upper bound by distance}, the area of region 
$\{X \in \Omega^\ast: |DX| \lesssim 2^{-n/2} r^{1/2} M^{1/2} \phi_L^{1/2} \}$ 
is dominated by $2^{-n/2} wr^{1/2} M^{1/2} \phi_L^{1/2}$ (up to a constant multiple). Furthermore, given such a point $X$, we may apply Lemma \ref{area upper bound by distance} again and obtain that the area of the region $\{Y\in \Omega^\ast: |XY|< 2L\}$ is dominated by $wL$. Finally, given a pair $(X,Y)$ as above, the area of the region $\{Z\in \Omega^\ast: \angle ZXY \in \Phi_n^L \}$ is dominated by $wr (2^{-n}\phi_L)$, thanks to Corollary \ref{area upper bound by sine}. A product of the three upper bounds above gives the upper bound of $|\tilde{\Omega}|$. In summary we always have 
\[
 \left|\Omega_{L,n}^{1,1} \right| \lesssim 2^{-n/2}  w r^{1/2} M^{1/2} \phi_L^{1/2} \cdot wL \cdot 2^{-n} wr \phi_L  \lesssim 2^{-3n/2} w^3 r^{3/2} M^{1/2} L \phi_L^{3/2}.
\]
Thus we have (In this case $\phi_L \lesssim L/r$)
\[
 \int_{\Omega_{L,n}^{1,1}} \frac{1}{|\triangle XYZ|} dX dY dZ \lesssim 2^{-n/2} w^3 r^{3/2} M^{1/2} L^{-1} \phi_L^{1/2} \lesssim 2^{-n/2} w^3 r M^{1/2} L^{-1/2}.
\]
If $L \geq \sqrt{wr}$, the case $k=2$ can be dealt with in the same way. We observe 
\[
 \Omega_{L,n}^{1,2} \subset  \left\{(X,Y,Z)\in (\Omega^\ast)^3: |DX| \lesssim M, \sin \angle DXZ \lesssim 2^{-n/2} r^{-1/2} M^{1/2} \phi_L^{1/2}, \angle ZXY \in \Phi_n^L\right\}
\]
We first choose an $X$ with $|DX| \lesssim M$, then determine the region containing all possible $Z$'s by the angle $\angle DXZ$, and finally determine the region of $Y$ by the angle $\angle ZXY$. This gives an upper bound 
\begin{align*}
 \int_{\Omega_{L,n}^{1,2}} \frac{1}{|\triangle XYZ|} dX dY dZ & \lesssim \frac{(wM)\cdot  (wr\cdot 2^{-n/2} r^{-1/2} M^{1/2} \phi_L^{1/2}) \cdot (wr \cdot 2^{-n} \phi_L) }{L^2\cdot  2^{-n} \phi_L}\\
  &\lesssim 2^{-n/2} w^3 r M^{3/2} L^{-3/2}. 
\end{align*}
We may deal with the case $M \ll L< \sqrt{wr}$ in exactly the same way by using Remark \ref{trivial area}, Lemma \ref{area upper bound by distance} and $\phi_L \lesssim w/L$. The upper bounds are given by
\begin{align*}
 \int_{\Omega_{L,n}^{1,1}} \frac{1}{|\triangle XYZ|} dX dY dZ & \lesssim \frac{(w\cdot 2^{-n/2} r^{1/2} M^{1/2} \phi_L^{1/2})\cdot (wL) \cdot (L^2 \cdot 2^{-n} \phi_L)}{L^2 \cdot 2^{-n} \phi_L}\\
 & \lesssim 2^{-n/2} w^{5/2} r^{1/2} M^{1/2} L^{1/2};\\
 \int_{\Omega_{L,n}^{1,2}} \frac{1}{|\triangle XYZ|} dX dY dZ & \lesssim \frac{(wM) \cdot (L^2 \cdot 2^{-n/2} r^{-1/2} M^{1/2} \phi_L^{1/2})\cdot (L^2 \cdot 2^{-n} \phi_L)}{L^2 \cdot 2^{-n} \phi_L} \\
 & \lesssim 2^{-n/2} w^{3/2} r^{-1/2} M^{3/2} L^{3/2}. 
\end{align*}
We may combine all the upper bounds above and conclude 
\[
 \sum_{L\gg M, n\geq 0} \int_{\Omega_{L,n}^1} \frac{1}{|\triangle XYZ|} dX dY dZ  \lesssim w^3 r. 
\]

\begin{figure}[h]
 \centering
 \includegraphics[scale=1.25]{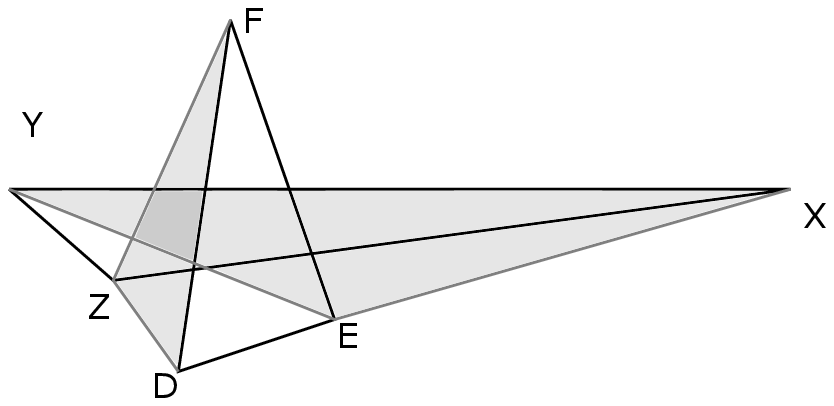}
 \caption{Large size, Type II reciprocal triangles} \label{figure large_type2}
\end{figure}

\paragraph{Type II} We may further write $\Omega_{L,n}^2 = \Omega_{L,n}^{2a} \cup \Omega_{L,n}^{2b}$ with
\begin{align*}
 \Omega_{L,n}^{2a} & = \left\{(X,Y,Z) \in \Omega_{L,n}^2: |ZD|, |ZE|, |ZF| \lesssim M\right\}; \\
 \Omega_{L,n}^{2b} & = \left\{(X,Y,Z) \in \Omega_{L,n}^2: |YD|, |YE|, |YF| \lesssim M\right\}.
\end{align*}
These two cases can be dealt with in exactly the same way. Let us consider the Type IIa reciprocal triangles, for instance. In this case
\[
  |XY|, |XZ|, |XD|, |XE|, |XF| \simeq L.
\]
By our reciprocal assumption, we have (see figure \ref{figure large_type2})
\[
 |\triangle FDZ| \cdot |\triangle EXY| \lesssim |\triangle DEF| \cdot |\triangle XYZ|. 
\]
That is
\[
 (|DZ|\cdot |DF| \sin \angle FDZ)(|XE|\cdot |XY|\sin \angle EXY) \lesssim (|DF|\cdot |FE|\sin \angle DFE) (|XY|\cdot |XZ|\sin \angle YXZ).
\]
Canceling $|DF|$, $|XY|$ and plugging $|FE|\simeq M$, $|XZ|, |XE| \simeq L$ in, we have 
\[
 |DZ| (\sin \angle FDZ)(\sin \angle EXY) \lesssim M (\sin \angle DFE) \sin \angle YXZ \lesssim M \phi_M \sin \angle YXZ.
\]
Following the same argument as in the Type I case, we may write $\Omega_{L,n}^{2a} = \Omega_{L,n}^{2a,1}\cup \Omega_{L,n}^{2a,2}\cup \Omega_{L,n}^{2a,3}$. 
Here we define 
\begin{align*}
 \Omega_{L,n}^{2a,1} & = \left\{(X,Y,Z) \in \Omega_{L,n}^{2a}: |DZ| \lesssim 2^{-n/3} r^{2/3} M^{1/3} \phi_M^{1/3} \phi_L^{1/3} \right\};\\
 \Omega_{L,n}^{2a,2} & = \left\{(X,Y,Z) \in \Omega_{L,n}^{2a}: \sin \angle FDZ \lesssim 2^{-n/3} r^{-1/3} M^{1/3} \phi_M^{1/3} \phi_L^{1/3} \right\};\\
\Omega_{L,n}^{2a,3} & = \left\{(X,Y,Z) \in \Omega_{L,n}^{2a}: \sin \angle EXY \lesssim 2^{-n/3} r^{-1/3} M^{1/3} \phi_M^{1/3} \phi_L^{1/3} \right\}.
\end{align*}
We then give upper bounds of the integrals below as in the Type I case: If $L \geq \sqrt{wr}$, then 
\begin{align*}
 \int_{\Omega_{L,n}^{2a,k}} \frac{1}{|\triangle XYZ|} dX dY dZ & \lesssim \frac{(w\cdot 2^{-n/3}r^{2/3}M^{1/3} \phi_M^{1/3} \phi_L^{1/3})\cdot (Lw) \cdot (wr \cdot 2^{-n} \phi_L)}{L^2 \cdot 2^{-n}\phi_L}\\
 & \lesssim 2^{-n/3} w^3 r^{4/3} L^{-2/3} (M\phi_M)^{1/3}. 
\end{align*}
On the other hand, if $M \ll  L < \sqrt{wr}$, then 
\begin{align*}
  \int_{\Omega_{L,n}^{2a,k}} \frac{1}{|\triangle XYZ|} dX dY dZ & \lesssim \frac{(w\cdot 2^{-n/3} r^{2/3} M^{1/3} \phi_M^{1/3} \phi_L^{1/3})\cdot (Lw) \cdot (L^2 \cdot 2^{-n} \phi_L)}{L^2 \cdot 2^{-n} \phi_L }\\
 & \lesssim 2^{-n/3} w^{8/3} r^{2/3} L^{2/3}.
\end{align*}
Finally we recall $M \phi_M \lesssim M^2/r$ if $M \geq \sqrt{wr}$ and $M \phi_M \lesssim w$ if $M \leq \sqrt{wr}$, then take a sum for all $L \gg M$ and $n\geq 0$. 
\[
 \sum_{L \gg M, n\geq 0} \int_{\Omega_{L,n}^{2a}} \frac{1}{|\triangle XYZ|} dX dY dZ \lesssim w^3r. 
\]
A similar inequality holds for Type IIb reciprocal triangles. 

\paragraph{Summary} We may combine Type I and II cases and obtain that for any given $\triangle DEF$, we have 
\[
 \sum_{L \gg M, n\geq 0} \int_{\Omega_{L,n}} \frac{1}{|\triangle XYZ|} dX dY dZ \lesssim w^3r. 
\]
Please note that the implicit constant in the inequality is an absolute constant, i.e. independent of $\triangle DEF$.  
\subsection{Small sizes}
We assume the size $L$ of $\triangle XYZ$ is much smaller than that of $\triangle DEF$, i.e. $L\ll M$. Again we consider Type I and II reciprocal triangles separately. We define  
\begin{align*}
 \Omega_{L,n}^1  & =  \left\{(X,Y,Z)\in \Omega_{L,n}: \triangle XYZ \; \hbox{and} \; \triangle DEF \; \hbox{are Type I reciprocal} \right\};\\
 \Omega_{L,n}^2  & =  \left\{(X,Y,Z)\in \Omega_{L,n}: \triangle XYZ \; \hbox{and} \; \triangle DEF \; \hbox{are Type II reciprocal} \right\}.
\end{align*}

\begin{figure}[h]
 \centering
 \includegraphics[scale=1.25]{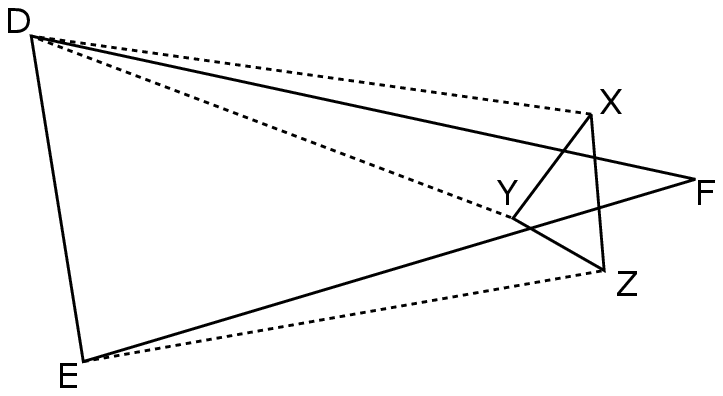}
 \caption{Small size, Type I reciprocal triangles} \label{figure small_type1}
\end{figure}

\paragraph{Type I} By our reciprocal assumption we always have (please see figure \ref{figure small_type1})
\[
 |\triangle XDY| \cdot |\triangle ZEF| \lesssim |\triangle XYZ| \cdot |\triangle DEF|.
\]
Thus
\[
 \left(|DY| \cdot |DX| \sin \angle XDY \right) \left(|EZ| \cdot |EF|\sin \angle ZEF\right) \lesssim \left(L^2 \sin \angle YXZ \right) |\triangle DEF|
\]
Our assumption implies $|DX|, |DY|, |EZ|, |EF| \simeq M$. Thus if $(X,Y,Z) \in \Omega_{L,n}^1$, we have 
\[
 (\sin \angle XDY)(\sin \angle ZEF) \lesssim 2^{-n} M^{-4} L^2 \phi_L  |\triangle DEF|.
\]
Thus we have $\Omega_{L,n}^1 = \Omega_{L,n}^{1,1} \cup \Omega_{L,n}^{1,2}$ with
\begin{align*}
 \Omega_{L, n}^{1,1}  & =  \left\{(X,Y,Z)\in \Omega_{L,n}^1: \sin \angle XDY \lesssim 2^{-n/2} M^{-2} L \phi_L^{1/2} |\triangle DEF|^{1/2}\right\}; \\
 \Omega_{L,n}^{1,2} & = \left\{(X,Y,Z)\in \Omega_{L,n}^1: \sin \angle ZEF \lesssim 2^{-n/2} M^{-2} L \phi_L^{1/2} |\triangle DEF|^{1/2}\right\}.
\end{align*}
If $M \gg L \geq \sqrt{rw}$, then we have (please note that $|XF| \lesssim L$, $\phi_L \lesssim L/r$ and $|\triangle DEF| \lesssim M^3/r$)
\begin{align*}
 \int_{\Omega_{L,n}^{1,k}} \frac{1}{|\triangle XYZ|} dX dY dZ &\lesssim \frac{(Lw) \cdot (wr\cdot 2^{-n/2} M^{-2} L \phi_L^{1/2} |\triangle DEF|^{1/2}) \cdot (wr \cdot 2^{-n} \phi_L)}{L^2 \phi_L 2^{-n}}\\
 & \lesssim \frac{2^{-n/2}w^3 r L^{1/2}}{M^{1/2}}. 
\end{align*}
If $L < \sqrt{rw} \leq M$, then we have ($\phi_L \lesssim w/L$)
\begin{align*}
 \int_{\Omega_{L,n}^{1,k}} \frac{1}{|\triangle XYZ|} dX dY dZ & \lesssim \frac{(Lw) \cdot (wr\cdot 2^{-n/2} M^{-2} L \phi_L^{1/2} |\triangle DEF|^{1/2}) \cdot (L^2 \cdot 2^{-n} \phi_L)}{L^2 \cdot 2^{-n}\phi_L} \\
 & \lesssim \frac{2^{-n/2}w^{5/2} r^{1/2} L^{3/2}}{M^{1/2}}. 
\end{align*}
Finally, if $L \ll M \leq \sqrt{wr}$, then we have ($\phi_L \lesssim w/L$; $|ZE|, |DY| \lesssim M$ and $|\triangle DEF| \lesssim Mw$)
\begin{align*}
 \int_{\Omega_{L,n}^{1,k}} \frac{1}{|\triangle XYZ|} dX dY dZ & \lesssim \frac{(Lw) \cdot (M^2\cdot 2^{-n/2} M^{-2} L \phi_L^{1/2} |\triangle DEF|^{1/2}) \cdot (L^2 \cdot 2^{-n} \phi_L)}{L^2 \cdot 2^{-n} \phi_L}\\ 
 & \lesssim 2^{-n/2} w^2 M^{1/2} L^{3/2}. 
\end{align*}
Collecting the upper bounds above and taking a sum, we always have 
\[
 \sum_{L \ll M, n \geq 0} \int_{\Omega_{L,n}^1} \frac{1}{|\triangle XYZ|} dX dY dZ \lesssim w^3 r. 
\]

\begin{figure}[h]
 \centering
 \includegraphics[scale=1.25]{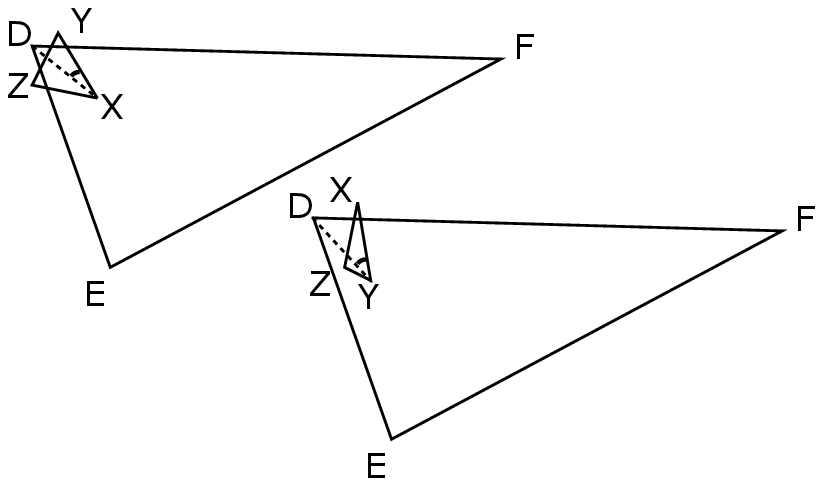}
 \caption{Small size, Type II reciprocal triangles} \label{figure small_type2}
\end{figure}

\paragraph{Type II} Now we consider small, type II reciprocal triangles of a given triangle $\triangle DEF$. This is the most difficult case. Let $\triangle XYZ$ of size $L$ be a Type II reciprocal triangle of $DEF$. Let us first give an upper bound of the integral
\[
 \int_{\Omega_{L,n}^2} \frac{1}{|\triangle XYZ|} dX dY dZ
\]
for given $L\ll M$, $n\geq 0$. Without loss of generality, let us assume \footnote{Strictly speaking, we need to consider four different cases. The argument given here only takes care of one from the four parts of $\Omega_{L,n}^2$. However, all these four cases can be dealt with in exactly the same way.}
\begin{align*}
 &|DX|, |DY|, |DZ| \lesssim L;& &|\triangle ZEF| \gtrsim |\triangle DEF|.&
\end{align*}
Thus by reciprocal assumption we immediately have 
\[
 |\triangle DXY| \lesssim |\triangle XYZ|. 
\]
Since $|XY| \simeq L$ and $|DX|, |DY|\lesssim L$, at least one of the following holds (see figure \ref{figure small_type2})
\begin{itemize}
 \item $|DX| \simeq L$. By comparing the area of $\triangle DXY$ with that of $\triangle XYZ$ we have
 \[ 
  |DX|\cdot |XY| \sin \angle DXY \lesssim |XY|\cdot |XZ| \sin \angle YXZ \; \Rightarrow \; \sin \angle DXY \lesssim \sin \angle YXZ.
 \]
 \item $|DY| \simeq L$. By considering the area of $\triangle DXY$ we have $\sin \angle DYX \lesssim \sin \angle YXZ$. 
\end{itemize}
Thus the region $\Omega_{L,n}^2$ is the union of two parts:
\begin{align*}
 \Omega_{L, n}^{2,1}  & =  \left\{(X,Y,Z)\in \Omega_{L,n}^2: \sin \angle DXY \lesssim 2^{-n} \phi_L \right\} \\
 \Omega_{L,n}^{2,2} & = \left\{(X,Y,Z)\in \Omega_{L,n}^2: \sin \angle DYX \lesssim 2^{-n} \phi_L \right\}.
\end{align*}
If $L\geq \sqrt{wr}$, we may find an upper bound of the integrals ($k=1,2$, $\phi_L \lesssim L/r$)
\begin{align}
 \int_{\Omega_{L,n}^{2,k}} \frac{1}{|\triangle XYZ|} dX dY dZ \lesssim \frac{(Lw) \cdot (wr\cdot 2^{-n} \phi_L) \cdot (wr \cdot 2^{-n} \phi_L)}{L^2 \cdot 2^{-n} \phi_L}  \lesssim 2^{-n} w^3 r^2 L^{-1} \phi_L \lesssim 2^{-n} w^3 r. \label{upper bound single L}
\end{align}
Similarly if $L < \sqrt{wr}$, then we have ($\phi_L \lesssim w/L$)
\begin{align} \label{upper bound single L2}
 \int_{\Omega_{L,n}^{2,k}} \frac{1}{|\triangle XYZ|} dX dY dZ \lesssim \frac{(Lw) \cdot (L^2 \cdot 2^{-n} \phi_L)\cdot (L^2\cdot 2^{-n} \phi_L )}{L^2 \cdot 2^{-n}\phi_L} \lesssim 2^{-n} w^2 L^2. 
\end{align}
We may collect the upper bounds above and obtain
\[
 \sum_{L \ll M, L \leq 32\sqrt{wr}} \sum_{n\geq 0}  \int_{\Omega_{L,n}^{2}} \frac{1}{|\triangle XYZ|} dX dY dZ \lesssim w^3 r;
\]
and 
\begin{align*}
  \sum_{32\sqrt{wr} < L \ll M} \left(\sum_{n > \log_2 \frac{\phi_L r^{1/2}}{8w^{1/2}}} \int_{\Omega_{L,n}^{2}} \frac{dX dY dZ}{|\triangle XYZ|} \right)
  & \lesssim \sum_{32\sqrt{wr} < L \ll M} \left(\sum_{n > \log_2 \frac{\phi_L r^{1/2}}{8w^{1/2}}} 2^{-n} w^3 r^2 L^{-1} \phi_L \right) \\
  & \lesssim \sum_{32\sqrt{wr} < L \ll M} w^{7/2} r^{3/2} L^{-1} \\
  & \lesssim w^3 r. 
\end{align*}
Thus it suffices to consider $(X,Y,Z) \in \Omega_{L,n}^2$ with $32\sqrt{wr} < L \ll M$ and $n \leq \log_2 (\phi_L r^{1/2} /8w^{1/2})$. We apply Lemma \ref{length of edge} and obtain 
\begin{align}
 |YZ| & \geq 2r \sin \angle YXZ - 2\sqrt{wr} - 2w \geq 2 r \phi_L 2^{-n-1} - 4\sqrt{wr} \geq 4\sqrt{wr}; \label{lower bound of YZ} \\
 |YZ| & \leq 2r\sin \angle YXZ + 2\sqrt{wr} \leq 2r \phi_L 2^{-n} + 2\sqrt{wr} \leq 3r \phi_L 2^{-n}. \label{upper bound of YZ by angle}
\end{align}
Next we first prove
\begin{lemma} \label{reclassification}
 Let $(X,Y,Z) \in \Omega_{L,n}^2$ with $32\sqrt{wr} < L\ll M$. In addition, we assume $|YZ| \geq 4\sqrt{wr}$. Then there exists an ansolute constant $c_1 >0$ so that at least one of the following holds
 \begin{itemize}
  \item[(a)] $c_1 |DE| \leq L \leq 8|DE|$; 
  \item[(b)] $L > 8|DE|$ and $\sin \angle EYX \lesssim 2^{-n} r^{-1} \max\{|DE|, \sqrt{wr}\}$; 
  \item[(c)] $L > 8|DE|$ and 
  \begin{align*}
   \min \{|DY|, |DZ|, |EY|, |EZ|\} \lesssim \min\{|YZ|, \max\{|DE|, \sqrt{wr}\}\}.
  \end{align*}
 \end{itemize}
\end{lemma}
\begin{proof}
The proof consists of three steps. 
\paragraph{Step 1}  We first show that $|DE| \lesssim L$. Without loss of generality we assume $|DX|, |DY|, |DZ| \lesssim L$ and $|XY|\geq L$. If $|DE| \gg L$, then we would have 
\[
 |EX|, |EY|, |EZ| \simeq |DE| \gg L > 32\sqrt{wr}. 
\]
Since $|XY| \geq L$, we have either $|DX| \geq  L/2$ or $|DY| \geq  L/2$. We consider these two cases separately. If $|DX| \geq L/2$, then our reciprocal assumption implies 
\[
 |\triangle DEX| \cdot |\triangle YZF| \lesssim |\triangle DEF| \cdot |\triangle XYZ|
\]
According to Corollary \ref{big triangle in omega}, the inequality above implies 
\[
 \frac{|DE| \cdot |EX| \cdot |DX|}{r}\cdot \frac{|YZ| \cdot |YF|\cdot |ZF|}{r} \lesssim \frac{|DE|\cdot |EF|\cdot |DF|}{r} \cdot \frac{|XY|\cdot |XZ|\cdot |YZ|}{r}
\]
We cancel $|YZ|, |DE|$, recall the facts
\[
|EX| \simeq |DE|,\qquad |DX|, |XY|, |XZ|\simeq L,\qquad |YF|, |ZF|, |EF|, |DF|\simeq M,
\]
and obtain $|DE|\lesssim L$. This is a contradiction. On the other hand, if $|DY| \geq L/2$, then we may follow a similar argument as above by considering $\triangle DEY, \triangle XZF$, and obtain
\[
 \frac{|DE| \cdot |EY| \cdot |DY|}{r}\cdot \frac{|XZ| \cdot |XF|\cdot |ZF|}{r} \lesssim \frac{|DE|\cdot |EF|\cdot |DF|}{r} \cdot \frac{|XY|\cdot |XZ|\cdot |YZ|}{r}
\]
This gives $|DE| \simeq |EY| \lesssim |YZ| \lesssim L$. Again this is a contradiction. As a result we obtain $|DE| \lesssim L$. It immediately follows that
\[
 |DX|, |DY|, |DZ|, |EX|, |EY|, |EZ| \lesssim L.
\]
Please refer to figure \ref{figure small_add1} for an illustration of the proof. Our remaining task is to show that if $|DE| < L/8$, then either (b) or (c) holds. 
 \begin{figure}[h]
 \centering
 \includegraphics[scale=1.25]{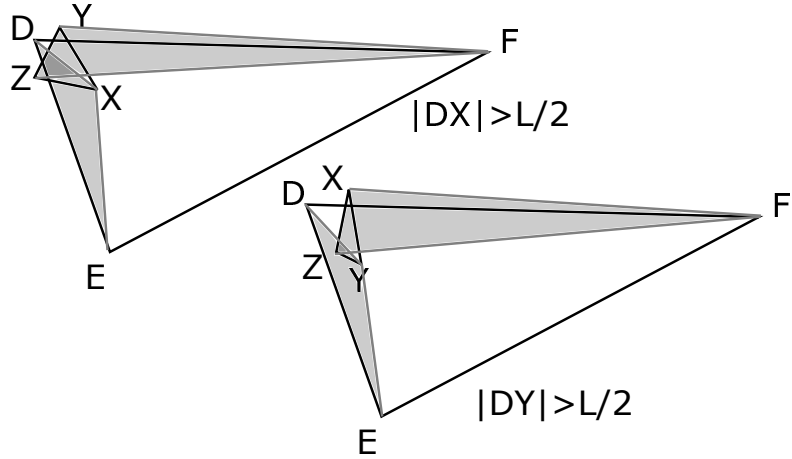}
 \caption{Large reciprocal triangles} \label{figure small_add1}
\end{figure}
\paragraph{Step 2} Now we assume $|DE| < L/8$, there are two cases: $D,E$ are either both close to the point $X$ or both far away from $X$. In this step we assume $|DX|, |XE| \leq L/4$. Since $|XY|, |XZ| > L/2$ we also have
\[
 |DY|, |EY|, |DZ|, |EZ| \geq L/4 > 8 \sqrt{wr}.
\]
We consider the triangles $\triangle EXY$ and $\triangle DZF$. The reciprocal assumption immediately gives
\[
 |\triangle EXY| \cdot |\triangle DZF| \lesssim |\triangle DEF| \cdot |\triangle XYZ|.
\]
Thus we may apply Corollary \ref{semi-big triangle in omega} and obtain
\begin{align*}
 (L^2 \sin \angle EYX) \cdot \frac{L M^2}{r} &\lesssim \frac{M^2 \max\{|DE|, \sqrt{wr}\}}{r} \cdot (L^2 \cdot 2^{-n}\phi_L) \\
 \Rightarrow \quad \sin \angle EYX & \lesssim 2^{-n} r^{-1} \max\{|DE|, \sqrt{wr}\}.
\end{align*}
In other words, (b) holds. Please see the upper half of figure \ref{figure small_add2}. 
 \begin{figure}[h]
 \centering
 \includegraphics[scale=1.25]{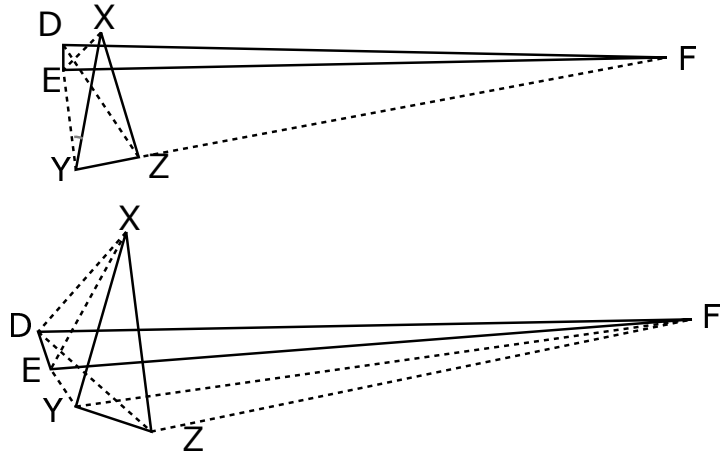}
 \caption{Type II reciprocal triangles of a narrow triangle} \label{figure small_add2}
\end{figure}
\paragraph{Step 3} Finally we assume $|DE| < L/8$ and $\max\{|DX|, |EX|\} > L/4$. This implies that $|DX|, |EX|\geq L/8 > 4\sqrt{wr}$. If we have 
\[
 \min\{|DY|, |DZ|, |EY|, |EZ|\} \leq 4\sqrt{wr},
\]
then our assumption on $|YZ|$ automatically guarantees (c) holds. Therefore we may additionally assume 
\[
 \min\{|DY|, |DZ|, |EY|, |EZ|\} > 4\sqrt{wr}.
\]
Without loss of generality, we assume 
\[
 |DZ| = \max\{|DY|, |DZ|, |EY|, |EZ|\}.
\]
We have
\begin{align*}
 &2|DZ| \geq |DZ|+|DY| \geq |YZ|;& &2|DZ| \geq |DZ|+|EZ| \geq |DE|.&
\end{align*}
Thus $|DZ| \geq |YZ|/2$, $|DZ| \geq \max\{|DE|, \sqrt{wr}\}/2$. We next apply Corollary \ref{big triangle in omega} and obtain
\begin{align*}
 &|\triangle XDZ| \gtrsim |\triangle XYZ|;& &|\triangle FDZ| \gtrsim |\triangle DEF|.&
\end{align*}
The reciprocal assumption then gives
\begin{align*} 
 &|\triangle EFY| \lesssim |\triangle DEF|;& &|\triangle EXY| \lesssim |\triangle XYZ|.&
\end{align*}
We then apply Corollary \ref{big triangle in omega} again and conclude (please refer to lower half of figure \ref{figure small_add2})
\begin{align*}
 &|EY| \lesssim \max\{|DE|, \sqrt{wr}\};& &|EY| \lesssim |YZ|.&
\end{align*}
Thus (c) holds. 
\end{proof}
\paragraph{Completion of type II case} First of all, we recall that it suffices to consider $(X,Y,Z) \in \Omega_{L,n}^2$ with $32\sqrt{wr} < L \ll M$ and $n\leq \log_2 (\phi_L r^{1/2} /8w^{1/2})$. According to Lemma \ref{reclassification}, the set $\Omega_{L,n}^2$ of this kind is empty unless $L \geq c_1 |DE|$. Thus we may further assume $L \geq c_1 |DE|$. We recall the upper bounds given in \eqref{upper bound single L}, \eqref{upper bound single L2} and obtain
\[
 \sum_{c_1 |DE| \leq L \leq 8|DE|} \sum_{n\geq 0} \int_{\Omega_{L,n}^2} \frac{1}{|\triangle XYZ|} dX dY dZ \lesssim w^3 r.
\]
Therefore we only need to deal with $\Omega_{L,n}^2$ with $\max\{32\sqrt{wr}, 8|DE|\} < L \ll M$ and $n\leq \log_2 (\phi_L r^{1/2} /8w^{1/2}) $. For convenience we use the notation $K = \max\{|DE|, \sqrt{wr}\}$. We recall \eqref{lower bound of YZ}, \eqref{upper bound of YZ by angle} and obtain that $(X,Y,Z) \in \Omega_{L,n}^2$ must satisfy $|YZ| \geq 4\sqrt{wr}$ and
\[
 \min\{|YZ|, \max\{|DE|, \sqrt{wr}\}\} \leq |YZ|^{1/2} K^{1/2} \lesssim 2^{-n/2} r^{1/2} \phi_L^{1/2} K^{1/2}.
\]
According to Lemma \ref{reclassification}, we may write 
\[
 \Omega_{L,n}^2 = \bigcup_{k=1}^5 \Omega_{L,n}^{2,k}.
\]
Here we define 
\begin{align*}
 \Omega_{L,n}^{2,1} & = \{(X,Y,Z)\in \Omega_{L,n}^2: \sin \angle EYX \lesssim 2^{-n} r^{-1} K\};\\
 \Omega_{L,n}^{2,2} & = \{(X,Y,Z)\in \Omega_{L,n}^2: |DY| \lesssim 2^{-n/2} r^{1/2} \phi_L^{1/2} K^{1/2}\};\\
 \Omega_{L,n}^{2,3} & = \{(X,Y,Z)\in \Omega_{L,n}^2: |DZ| \lesssim 2^{-n/2} r^{1/2} \phi_L^{1/2} K^{1/2}\};\\
 \Omega_{L,n}^{2,4} & = \{(X,Y,Z)\in \Omega_{L,n}^2: |EY| \lesssim 2^{-n/2} r^{1/2} \phi_L^{1/2} K^{1/2}\};\\
 \Omega_{L,n}^{2,5} & = \{(X,Y,Z)\in \Omega_{L,n}^2: |EZ| \lesssim 2^{-n/2} r^{1/2} \phi_L^{1/2} K^{1/2}\}.
\end{align*}
We then apply Lemma \ref{area upper bound by distance}, Corollary \ref{area upper bound by sine} and obtain ($k=2,3,4,5$)
\begin{align*}
 \int_{\Omega_{L,n}^{2,1}} \frac{1}{|\triangle XYZ|} dX dY dZ & \lesssim \frac{(wL) \cdot (wr \cdot 2^{-n} r^{-1} K) \cdot (wr \cdot 2^{-n}\phi_L)}{L^2 \cdot 2^{-n} \phi_L} \lesssim 2^{-n} w^3 r L^{-1} K;\\
 \int_{\Omega_{L,n}^{2,k}} \frac{1}{|\triangle XYZ|} dX dY dZ & \lesssim \frac{(w \cdot 2^{-n/2} r^{1/2} \phi_L^{1/2} K^{1/2})\cdot (wL) \cdot (wr \cdot 2^{-n} \phi_L)}{L^2 \cdot 2^{-n}\phi_L}\\
 & \lesssim 2^{-n/2} w^3 r K^{1/2} L^{-1/2}.
\end{align*}
Thus 
\[
 \sum_{\max\{32\sqrt{wr}, 8|DE|\} < L \ll M} \left(\sum_{n \leq \log_2 \frac{\phi_L r^{1/2}}{8w^{1/2}}} \int_{\Omega_{L,n}^{2}} \frac{1}{|\triangle XYZ|} dX dY dZ \right) \lesssim w^3 r. 
\]
In summary we have 
\[
 \sum_{L\ll M, n\geq 0} \int_{\Omega_{L,n}^2} \frac{1}{|\triangle XYZ|} dX dY dZ \lesssim w^3 r.
\]
\paragraph{Summary} We may combine Type I and II cases and obtain that 
\[
 \sum_{L \ll M, n\geq 0} \int_{\Omega_{L,n}} \frac{1}{|\triangle XYZ|} dX dY dZ \lesssim w^3 r.
\]
\subsection{Comparable Sizes}

\begin{figure}[h]
 \centering
 \includegraphics[scale=1.25]{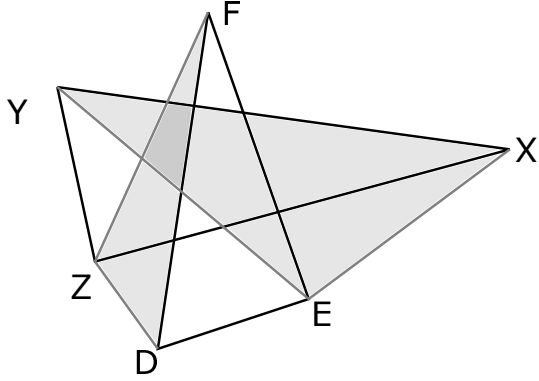}
 \caption{Comparable size reciprocal triangles} \label{figure equal_uni}
\end{figure}

Finally let us the consider the case when $\triangle XYZ$ and $\triangle DEF$ are about of the same size, i.e. $L\simeq M$. This eliminate the need to take a sum in $L$. In this subsection we prove that if $L \simeq M$, then
\[
 \sum_{n\geq 0} \int_{\Omega_{L,n}} \frac{1}{|\triangle XYZ|} dX dY dZ \lesssim w^3 r. 
\]
The argument is similar to the case $L \gg M$, Type II. Now we have less information on the relative location of two triangles available. Nevertheless, Corollary \ref{close to each other} guarantees that $d(\triangle XYZ, \triangle DEF) \lesssim L\simeq M$. By reciprocal assumption, we have (please refer to figure \ref{figure equal_uni}) 
\[
 |\triangle FDZ| \cdot |\triangle EXY| \lesssim |\triangle DEF| \cdot |\triangle XYZ|. 
\]
That is
\[
 (|DZ|\cdot |DF| \sin \angle FDZ)(|XE|\cdot |XY|\sin \angle EXY) \lesssim (|DF|\cdot |FE|\sin \angle DFE) (|XY|\cdot |XZ|\sin \angle YXZ).
\]
Canceling $|DF|$, $|XY|$ and plugging $|FE|, |XZ|\simeq M$ in, we have 
\begin{align*}
 |DZ| (\sin \angle FDZ)|XE|(\sin \angle EXY) & \lesssim M^2 (\sin \angle DFE) \sin \angle YXZ \\
 & \lesssim M^2 \phi_M \cdot 2^{-n} \phi_L\\
 & \lesssim 2^{-n} K_M^2.
\end{align*}
Here the notation $K_M$ represents
\[
 K_M = \left\{\begin{array}{ll} M^2/r, & \hbox{if}\; M\geq \sqrt{wr}; \\ w, & \hbox{if}\; M < \sqrt{wr}.\end{array}\right.
\]
Therefore we may write $\Omega_{L,n} = \Omega_{L,n}^{0,1} \cup \Omega_{L,n}^{0,2} \cup \Omega_{L,n}^{0,3} \cup \Omega_{L,n}^{0,4}$. Here 
\begin{align*}
 \Omega_{L,n}^{0,1} & = \left\{(X,Y,Z) \in \Omega_{L,n}: |DZ| \lesssim 2^{-n/4} r^{1/2} K_M^{1/2}\right\};\\
 \Omega_{L,n}^{0,2} & = \left\{(X,Y,Z) \in \Omega_{L,n}: \sin \angle FDZ \lesssim 2^{-n/4} r^{-1/2} K_M^{1/2}\right\};\\
 \Omega_{L,n}^{0,3} & = \left\{(X,Y,Z) \in \Omega_{L,n}: |XE| \lesssim 2^{-n/4} r^{1/2} K_M^{1/2}\right\}; \\
\Omega_{L,n}^{0,4} & = \left\{(X,Y,Z) \in \Omega_{L,n}: \sin \angle EXY \lesssim 2^{-n/4} r^{-1/2} K_M^{1/2}\right\}.
\end{align*}
This immediately gives the upper bounds: if $M \geq \sqrt{wr}$, then 
\begin{align*}
 \int_{\Omega_{L,n}^{0,k}} \frac{1}{|\triangle XYZ|} dX dY dZ \lesssim \frac{(w\cdot 2^{-n/4} r^{1/2} K_M^{1/2})\cdot (Lw) \cdot (wr\cdot 2^{-n}\phi_L )}{L^2 \cdot 2^{-n}\phi_L} \lesssim 2^{-n/4} w^3 r.
\end{align*}
If $M \leq \sqrt{wr}$, then 
\begin{align*}
 \int_{\Omega_{L,n}^{0,k}} \frac{1}{|\triangle XYZ|} dX dY dZ & \lesssim \frac{(w\cdot 2^{-n/4} r^{1/2} K_M^{1/2})\cdot (Lw) \cdot (L^2 \cdot 2^{-n}\phi_L)}{L^2 \cdot 2^{-n}\phi_L} \\
 &\lesssim 2^{-n/4} w^{5/2} r^{1/2} L \lesssim 2^{-n/4} w^3 r.
 \end{align*}
In either case we may take a sum and obtain that if $L \simeq M$, then
\[
 \sum_{n\geq 0} \int_{\Omega_{L,n}} \frac{1}{|\triangle XYZ|} dX dY dZ \lesssim w^3 r.
\]

\subsection{Summary} 
Collecting all cases discussed above, we prove that the inequality
\[
 \int_{\Sigma(DEF)\cap (\Omega^\ast)^3} \frac{1}{|\triangle XYZ|} dX dY dZ \lesssim w^3 r.
\]
holds for all $D,E,F\in \Omega^\ast$. The implicit constant here is an absolute constant. Thus we finish the proof of Proposition \ref{geometric inequality}. 

\section{Applications of Geometric Inequalities}
In this section we prove the main results given in Section \ref{sec:intro}. 
\subsection{Proof of Proposition \ref{main 1}} 
\paragraph{Part (a)} Let us temporally assume $G$ is supported in $[a,b]\times \mathbb{S}^2$. We apply Proposition \ref{geometric inequality} and obtain an upper bound of $C_{a,b,\delta}$ defined in Proposition \ref{transform to geometric problem}. (A disk of radius $r$ can be viewed as an annulus of outer radius $r$ and width $r$)
\[
 C_{a,b,\delta}(h) \lesssim \left\{\begin{array}{ll} (b/h)^6 \cdot w^3 r, & h\geq b; \\ (b/h)^6 \cdot r^4, & h \in (a-\delta, b).\end{array} \right.
\]
Here 
\begin{align*}
 r & = \frac{\sqrt{h^2-(a-\delta)^2}}{a-\delta};\\
 w & = \frac{\sqrt{h^2-(a-\delta)^2}}{a-\delta} - \frac{\sqrt{h^2-b^2}}{b}.
\end{align*}
We plug $r,w$ and obtain that if $h \geq b$, then (recall that $b/a \leq 2$ and $\delta>0$ is small)
\[
 C_{a,b,\delta}(h) \lesssim \frac{b^6 [1/(a-\delta)-1/b]^3}{h^2 (a-\delta)^2 [1/(a-\delta)-1/h]} \lesssim \frac{a^4 [1/(a-\delta)-1/b]^3}{h^2 [1/(a-\delta)-1/h]}.
\]
And if $h \in (a-\delta,b)$, then (in this case $a-\delta \simeq b \simeq h$)
\[
 C_{a,b,\delta}(h) \lesssim h^2[1/(a-\delta)-1/h]^2.
\]
This upper bound of $C_{a,b,\delta}(h)$ is an increasing function of $h$ in the interval $(a-\delta,b)$ and a decreasing function of $h$ in the interval $[b,+\infty)$. Thus we have 
\begin{align*}
 \sup_{h>\max\{a-\delta, R\}} C_{a,b,\delta}(h) & \lesssim \frac{a^4 [1/(a-\delta)-1/b]^3}{R^2 [1/(a-\delta)-1/R]}, & & R \geq b; \\
 \sup_{h>\max\{a-\delta, R\}} C_{a,b,\delta}(h) & \lesssim b^2 [1/(a-\delta) - 1/b]^2, & & R < b.
\end{align*}
We plug this upper bound in Proposition \ref{transform to geometric problem}, make $\delta\rightarrow 0^+$ and obtain 
\[
 \int_{|x|>R} |\mathbf{T} G (x)|^6 dx \lesssim \frac{a^4 (1/a-1/b)^3}{R^2 (1/a-1/R)} \|G\|_{L^2(\Rm\times \mathbb{S}^2)}^6 = \frac{(a/R)^2(1-a/b)^3}{1-a/R} \|G\|_{L^2(\Rm\times \mathbb{S}^2)}^6
\]
for all $L^2$ functions $G$ supported in $[a,b]\times \mathbb{S}^2$ and $R \geq b$. Similarly we may choose $R =0$, recall $1<b/a\leq 2$ and obtain
\[
 \int_{\Rm^3} |\mathbf{T} G (x)|^6 dx \lesssim (1-a/b)^2 \|G\|_{L^2(\Rm\times \mathbb{S}^2)}^6.
\]
By the identity
\[
 \mathbf{T} (G(-s,\omega))(x) = (\mathbf{T} G(s,\omega))(-x),
\]
The same inequalities as above also hold for $G$ supported in $[-b,-a]\times \mathbb{S}^2$. We then use the linearity of $\mathbf{T}$ to finish the proof.

\paragraph{Part (b)} Now let us assume $G\in L^2(\Rm\times \mathbb{S}^2)$ is supported in $[-b,b] \times \mathbb{S}^2$. We may break $G$ into pieces
\[
 G(s,\omega) = \sum_{k=0}^\infty G_k (s,\omega) 
\]
so that 
\[
  G_k(s,\omega) = \left\{\begin{array}{ll} G(s,\omega), & 2^{-k-1} b < |s| \leq 2^{-k} b; \\ 0, & \hbox{otherwise}. \end{array}\right. 
\]
It immediately gives a convergence in $L^6(\Rm^3)$: 
\[
 \mathbf{T} G =  \sum_{k=0}^\infty \mathbf{T} G_k = \sum_{k=0}^\infty \int_{\mathbb{S}^2} G_k (x\cdot \omega, \omega) d\omega. 
\]
We then apply the conclusion of part (a) on the radiation profiles $G_k$ and obtain
\begin{align*}
 \int_{|x|>R} \left|\mathbf{T} G_k (x)\right|^6 dx  \lesssim  \frac{\left(\frac{2^{-k-1}b}{R}\right)^2 \left(1 - \frac{2^{-k-1}b}{2^{-k} b}\right)^3}{1-\frac{2^{-1-k}b}{R}}  \|G_k\|_{L^2}^6  \lesssim \frac{2^{-2k} b^2}{R^2} \|G\|_{L^2}^6.
\end{align*}
Therefore 
\[
 \|\mathbf{T} G\|_{L^6(\{x: |x|>R\})} \leq \sum_{k=0}^\infty \|\mathbf{T} G_k\|_{L^6(\{x: |x|>R\})} \lesssim (b/R)^{1/3} \|G\|_{L^2}. 
\]
This finishes the proof of part (b).

\subsection{Proof of Corollary \ref{cor1}}
 
Since we always have $\|\mathbf{T} G\|_{L^6(\Rm^3)} \lesssim \|G\|_{L^2}$, we may assume $R > 2(b-a)$, without loss of generality. Let $c= \max\{|a|,|b|\}$ thus we have $\hbox{Supp}\, G \subset [-c,c]$. There are two cases
\begin{itemize}
 \item If $b-a \geq c/2$, then we apply part (b) of Proposition \ref{main 1} and obtain 
 \[
  \int_{|x|>R} |\mathbf{T} G(x)|^6 dx \lesssim (c/R)^2 \|G\|_{L^2}^6 \lesssim \left(\frac{b-a}{R}\right)^2 \|G\|_{L^2}^6. 
 \]
 \item If $b-a < c/2$, then we have either $2a > b > a > 0$ or $2b < a < b < 0$. Let us consider the first situation since the second one can be dealt with by symmetry. We apply Part (a) of Proposition \ref{main 1} and obtain 
 \begin{align*}
  \int_{|x|>R} |\mathbf{T} G(x)|^6 dx  & \lesssim (1-a/b)^2 \|G\|_{L^2}^6 \leq \left(\frac{b-a}{R}\right)^2 \|G\|_{L^2}^6, & &R \leq b;&\\
 \int_{|x|>R} |\mathbf{T} G(x)|^6 dx & \lesssim  \frac{(a/R)^2 (1-a/b)^3}{1-a/R} \|G\|_{L^2}^6 \lesssim \left(\frac{b-a}{R}\right)^2 \|G\|_{L^2}^6, & & R > b.&
 \end{align*}
 
\end{itemize}

\subsection{Proof of Proposition \ref{main 2}}

\paragraph{The $L^6$ estimates} Let $G$ be the radiation profile associated to the linear free wave $u$. By isometric property we have $E = 2\|G\|_{L^2}^2$. The non-radiative assumption implies that $G$ is supported in $[-r,r]\times \mathbb{S}^2$. We may define $G^{(t)}(s,\omega) = G(s+t, \omega)$ and rewrite 
 \[
  u(x,t) = \frac{1}{2\pi}\int_{\mathbb{S}^2} G(x\cdot \omega+t, \omega) d\omega = \frac{1}{2\pi}\int_{\mathbb{S}^2} G^{(t)} (x\cdot \omega, \omega) d\omega.
 \]
 We have 
 \begin{align*}
  &\hbox{Supp} \,G^{(t)} (s,\omega) \subseteq [-t-r,-t+r]\times \mathbb{S}^2;& &\|G^{(t)}\|_{L^2(\Rm\times \mathbb{S}^2)} = \|G\|_{L^2(\Rm\times \mathbb{S}^2)}.&  
 \end{align*}
We then apply Corollary \ref{cor1} and obtain
\[
 \|u(\cdot,t)\|_{L^6(\{x:|x|>R\})} \lesssim \left(2r/R\right)^{1/3} \|G^{(t)}\|_{L^2} \lesssim (r/R)^{1/3} \|G\|_{L^2}. 
\]
Therefore we have 
\[
 \|u\|_{L^\infty L^6(\Rm \times \{x:|x|>R\})} \lesssim (r/R)^{1/3} \|G\|_{L^2} \simeq (r/R)^{1/3} E^{1/2}. 
\]
The decay $\|u(\cdot,t)\|_{L^6(\Rm^3)} \lesssim (r/|t|)^{1/3} E^{1/2}$ immediately follows Proposition \ref{main 1}, part (a), as long as $|t|>3r$
\[
 \|u(\cdot,t)\|_{L^6(\Rm^3)} \lesssim \left(1-\frac{|t|-r}{|t|+r}\right)^{1/3} \|G^{(t)}\|_{L^2} \lesssim (r/|t|)^{1/3} E^{1/2}. 
\]
The case $|t|\leq 3r$ is trivial. 

\paragraph{The $L^p L^q$ estimates} We recall the Strichartz estimates given in Ginibre-Velo \cite{strichartz}: if $p_1, q_1 > 0$ satisfies $1/p_1 + 3/q_1 = 1/2$ and $6 \leq q_1 < +\infty$, then any finite-energy linear free wave $u$ satisfies 
\begin{equation} \label{regular Strichartz}
 \|u\|_{L^{q_1} L^{q_1}(\Rm \times \Rm^3)} \lesssim_{p_1,q_1} E^{1/2}. 
\end{equation} 
As a result, the $(r/R)^{\kappa}$ decay of the norm $\|u\|_{L^p L^q(\Rm \times \{x: |x|>R\})}$ follows an interpolation between the decay estimate 
\[
 \|u\|_{L^\infty L^6(\Rm \times \{x:|x|>R\})} \lesssim (r/R)^{1/3} E^{1/2}
\]
and the regular Strichartz estimate  \eqref{regular Strichartz} with $p_1 = 2^+$ and $q_1 = \infty^-$ in the whole space. Please note that the choice $(p_1,q_1) = (2, +\infty)$ is forbidden in the Strichartz estimates. 

\section{Two dimensional case and application on Radon transform}
 Since $\mathbf{T}$ is the adjoint operator of Radon transform $\mathcal{R}$, a corollary immediately follows Proposition \ref{main 1}: If $f$ is supported in the region $\{x: |x|>R\}$, then
 \[
  \|\mathcal{R} f\|_{L^2([-b,b]\times \mathbb{S}^2)} \lesssim (b/R)^{1/3} \|f\|_{L^{6/5}(\Rm^3)}, \qquad \forall \; b \in (0,R). 
 \]
We are also interested in its 2-dimensional analogue, since 2-dimensional Radon transform is more frequently used in some applications, for example, X-ray technology. Our 2-dimensional result is 
\begin{proposition}
 We consider the 2-dimensional Radon transform ($dL$ is the line measure of a straight line on $\Rm^2$)
\[
 \mathcal{R} f (s,\omega)= \int_{\omega \cdot x = s} f(x) dL(x), \qquad (s,\omega) \in \Rm \times \mathbb{S}^1; 
\]
and its adjoint 
\[
 \mathcal{R}^\ast G(x) = \int_{\mathbb{S}^1} G(x\cdot \omega, \omega) d\omega, \qquad x\in \Rm^2.
\]
Then we have 
\begin{itemize}
 \item[(a)] Assume that $b>a>0$ with $b/a\leq 2$. If $G$ is supported in $([-b,-a]\cup[a,b])\times \mathbb{S}^1$, then 
 \begin{align*}
  \int_{|x|>R} |\mathcal{R}^\ast G (x)|^4 dx &\lesssim \frac{(a/R)(1-a/b)}{(1-a/R)^{1/2}} \|G\|_{L^2(\Rm\times \mathbb{S}^1)}^4, \qquad \forall R>b; \\
  \int_{\Rm^2}  |\mathcal{R}^\ast G (x)|^4 dx &\lesssim (1-a/b)^{1/2} \|G\|_{L^2(\Rm\times \mathbb{S}^1)}^4.
 \end{align*}
 \item[(b)] Assume that $R \geq b > 0$. If $G$ is supported in $([-b,b])\times \mathbb{S}^1$, then
 \[
  \int_{|x|>R} |\mathcal{R}^\ast G (x)|^4 dx \lesssim \frac{b}{R}\|G\|_{L^2(\Rm\times \mathbb{S}^1)}^4.
 \] 
 \item[(c)] If $f(x) = 0$ for all $|x|<R$, then we have 
 \[
  \|\mathcal{R} f\|_{L^2([-b,b]\times \mathbb{S}^1)} \lesssim (b/R)^{1/4} \|f\|_{L^{4/3} (\Rm^2)}. 
 \]
\end{itemize}
\end{proposition}
\begin{proof}
 The general idea is exactly the same as in the 3-dimensional case. Part (b) follows part (a) and a decomposition of $G$. Part (c) follows a basic property of adjoint operators. Let us stretch the proof of part (a) only. Following the same argument as in Section \ref{sec: transformation}, we obtain ($\delta >0$ is sufficiently small)
\[
 \int_{|x|>R} |\mathbf{T} G (x)|^4 dx \leq \left(\sup_{h>\max\{R,a-\delta\}} C_{a,b,\delta}(h) \right)\|G\|_{L^2(\Rm\times \mathbb{S}^1)}^4.
\]
Here the constant $C_{a,b,\delta}(h)$ is defined by
\[
 C_{a,b,\delta}(h) = \frac{3b^2}{h^2}\sup_{x_1, x_2 \in \Omega_{\delta,h}^\ast} \int_{\Sigma(x_1,x_2)\cap \Omega_{\delta,h}^\ast} \frac{1}{|x_3-x_4|} dx_3 dx_4
\]
Here $\Omega_{\delta,h}^\ast$ is defined by
\begin{align*}
 \Omega_{\delta,h}^\ast & = \left(-\frac{\sqrt{h^2-(a-\delta)^2}}{a-\delta}, -\frac{\sqrt{h^2-b^2}}{b}, \right) \cup \left(\frac{\sqrt{h^2-b^2}}{b}, \frac{\sqrt{h^2-(a-\delta)^2}}{a-\delta}\right).& & h\geq b;\\
 \Omega_{\delta,h}^\ast & = \left(-\frac{\sqrt{h^2-(a-\delta)^2}}{a-\delta},\frac{\sqrt{h^2-(a-\delta)^2}}{a-\delta}\right), & & h\in (a-\delta,b);
\end{align*}
and $\Sigma(x_1,x_2)$ is the set consisting of reciprocal pairs of $(x_1,x_2)$. We call $(x_1,x_2)$ and $(x_3,x_4)$ are reciprocal pairs if and only if 
\[
 |x_1-x_2|\cdot |x_3-x_4| \geq \frac{1}{17} \max_{\{j_1,j_2, j_3,j_4\}=\{1,2,3,4\}} |x_{j_1}-x_{j_2}|\cdot |x_{j_3}-x_{j_4}|. 
\]
We claim that if $\Omega = (-r,-r+w) \cup (r-w,r)$ with $r\geq w>0$, then 
\begin{equation} \label{2d geometric inequality} 
 \sup_{x_1, x_2 \in \Omega} \int_{\Sigma(x_1,x_2)\cap \Omega^2} \frac{1}{|x_3-x_4|} dx_3 dx_4 \lesssim 
 w.
\end{equation}
We then plug this upper bound in the expression of $C_{a,b,\delta}(h)$, take the least upper bound, then make $\delta \rightarrow 0^+$ to finish the proof of part (a). Finally we need to verify \eqref{2d geometric inequality}. First of all, if $(x_1,x_2)$ and $(x_3,x_4)$ are reciprocal pairs, then we claim
\begin{equation} \label{necessary reciprocal}
 \min\{|x_1-x_3|, |x_1-x_4|, |x_2-x_3|, |x_2-x_4|\} \leq 34 \min\{|x_1-x_2|, |x_3-x_4|\}. 
\end{equation}
In fact, we may assume $|x_3-x_4|\geq |x_1-x_2|$ without loss of generality. The triangle inequality implies that we have either $|x_1-x_3| \geq |x_3-x_4|/2$ or $|x_1-x_4| \geq |x_3-x_4|/2$. By our reciprocal assumption we have either $|x_2-x_4|\leq 34 |x_1-x_2|$ or $|x_2-x_3|\leq 34|x_1-x_2|$. This verifies \eqref{necessary reciprocal}. Now we are ready to prove \eqref{2d geometric inequality}. Let us fix $x_1,x_2 \in \Omega$. We define $\Omega_L$ to be the set of all reciprocal pairs of size $L$:
\[
 \Omega_L = \{(x_3,x_4)\in \Sigma(x_1,x_2)\cap \Omega^2: L \leq |x_3-x_4|<2L\}, \qquad L\in \{r,r/2,r/4,\cdots\}. 
\]
We also assume that the size of $|x_1-x_2|$ is $M$, i.e. $M \in \{r,r/2,r/4,\cdots\}$ so that $M \leq |x_1-x_2|<2M$. There are two cases: Case 1, if $L \ll M\leq r$, then we must have either $L \leq w$ or $\Omega_L = \varnothing$. In addition \eqref{necessary reciprocal} implies that if $(x_3,x_4)\in \Omega_L$, then either $|x_4-x_1|, |x_3 - x_1| \lesssim L$ or $|x_4-x_2|, |x_3 - x_2|\lesssim L$ holds, thus $|\Omega_L| \lesssim L^2$. This implies 
\[
 \int_{\Omega_L} \frac{1}{|x_3-x_4|} dx_3 dx_4 \lesssim L^{-1} |\Omega_L|\lesssim L. 
\]
Case 2, if $L \gtrsim M$, then we have either $|x_3 - x_1| \lesssim M$ or $|x_4 - x_1|\lesssim M$ by \eqref{necessary reciprocal}, thus $|\Omega_L| \lesssim Mw$. (Please note that $|\Omega| = 2w$) Therefore we have 
\[
 \int_{\Omega_L} \frac{1}{|x_3-x_4|} dx_3 dx_4 \lesssim MwL^{-1}.
\]
In summary
\[
 \int_{\Sigma(x_1,x_2)\cap \Omega^2} \frac{dx_3 dx_4}{|x_3-x_4|} \leq \sum_{L} \int_{\Omega_L} \frac{dx_3 dx_4}{|x_3-x_4|} \lesssim \sum_{L\ll M, L\leq w} L + \sum_{L \gtrsim M} MwL^{-1}\lesssim w. 
\]
This finishes the proof. 
\end{proof}

\section*{Acknowledgement}
The authors are financially supported by National Natural Science Foundation of China Project 12071339.

\end{document}